\documentclass[12pt]{article}

\usepackage{subcaption}

\usepackage{float}

\usepackage{multirow}

\usepackage{tikz}

\usetikzlibrary{calc,backgrounds,arrows,matrix}

\usepackage{enumerate}

\usepackage{color}
\definecolor{lightblue}{rgb}{0,0.2,0.5}

\usepackage[colorlinks=true, urlcolor=blue,linkcolor=blue, citecolor=lightblue]{hyperref}

\usepackage[round,sort,comma,numbers]{natbib}

\bibpunct{\textcolor{lightblue}{(}}{\textcolor{lightblue}{)}}{,}{a}{}{;}

\usepackage{amssymb,amsmath}

\usepackage{graphicx}

\usepackage{graphics}
\usepackage{color}

\usepackage{tikz}
\usetikzlibrary{trees}
\usetikzlibrary{shapes.geometric} 
\usetikzlibrary{shapes,snakes}

\DeclareMathAlphabet{\eufrak}{U}{}{}{}
\SetMathAlphabet\eufrak{normal}{U}{euf}{m}{n}
\SetMathAlphabet\eufrak{bold}{U}{euf}{b}{n}

\allowdisplaybreaks

 \def\qu{{\mathord{\mathbb Z}}}

 \def\T{{\mathrm{{\rm T}}}}

 \def\inte{{\mathord{\mathbb R}}}

 \def\inte{{\mathord{\mathbb N}}}

 \def\sZZ{{\rm Z\kern-.45em{}Z}}

 \def\sQQ{{\kern 0.27em \vrule height1.45ex width0.03em depth0em
           \kern-0.30em \rm Q}}
 \def\qu{{\mathchoice
         {\sQQ}
         {\sQQ}
   {\kern 0.225em \vrule height1.05ex width0.025em depth0em \kern-0.25em \rm Q}
   {\kern 0.180em \vrule height0.78ex width0.020em depth0em \kern-0.20em \rm Q}
         }}
 \def\sGG{{\kern 0.27em \vrule height1.45ex width0.03em depth0em
           \kern-0.30em \rm G}}
 \def\gg{{\mathchoice
         {\sGG}
         {\sGG}
   {\kern 0.225em \vrule height1.05ex width0.025em depth0em \kern-0.25em \rm G}
   {\kern 0.180em \vrule height0.78ex width0.020em depth0em \kern-0.20em \rm G}
         }}

 \newtheorem{prop}{Proposition}[section]
 \newtheorem{lemma}[prop]{Lemma}
 \newtheorem{definition}[prop]{Definition}
 
 \newtheorem{theorem}[prop]{Theorem}
 \newtheorem{remark}[prop]{Remark}

\numberwithin{equation}{section}

 \def\P{{\mathord{\mathbb P}}}

 \newcounter{hyp}
\newenvironment{Proofy}{\removelastskip\par\medskip
\noindent{\em Proof of Theorem} \rm}{\penalty-20\null\hfill$\square$\par\medbreak}

\newenvironment{Proof}{\removelastskip\par\medskip \noindent{\em Proof.} \rm}{\penalty-20\null\hfill$\square$\par\medbreak}

\def\bprf{\begin{Proof}}
\def\nprf{\end{Proof}}
\def\bdes{\begin{description}}
\def\ndes{\end{description}}

\newtheorem{thm}{Theorem}[section]

\def\bdef{\begin{defn}}
\def\ndef{\end{defn}}
\def\bthm{\begin{thm}}
\def\nthm{\end{thm}}
\def\bprop{\begin{prop}}
\def\nprop{\end{prop}}
\def\brmk{\begin{remark}}
\def\nrmk{\end{remark}}
\def\bexa{\begin{exa}}
\def\nexa{\end{exa}}
\def\blem{\begin{lem}}
\def\nlem{\end{lem}}
\def\bcor{\begin{cor}}
\def\ncor{\end{cor}}
\def\bexe{\begin{exe}}
\def\nexe{\end{exe}}

\newcommand{\E}{\mathbb{E}}

\newcommand{\real}{\mathbb{R}}

\def\og{\leavevmode\raise.3ex
     \hbox{$\scriptscriptstyle\langle\!\langle$~}}
\def\fg{\leavevmode\raise.3ex
     \hbox{~$\!\scriptscriptstyle\,\rangle\!\rangle$}~}

\title{\Huge
Existence of solutions for nonlinear elliptic PDEs with fractional Laplacians on open balls
}

\author{
  Guillaume Penent\footnote{\href{mailto:PENE0001@e.ntu.edu.sg}{pene0001@e.ntu.edu.sg}}
  \qquad 
      Nicolas Privault\footnote{
\href{mailto:nprivault@ntu.edu.sg}{nprivault@ntu.edu.sg}
}
  \\
\small
Division of Mathematical Sciences
\\
\small
School of Physical and Mathematical Sciences
\\
\small
Nanyang Technological University
\\
\small
21 Nanyang Link, Singapore 637371
}

\allowdisplaybreaks

\usepackage{empheq}

\usepackage{bbm}

\makeatletter
\newcommand*\rel@kern[1]{\kern#1\dimexpr\macc@kerna}
\newcommand*\widebar[1]{
  \begingroup
  \def\mathaccent##1##2{
    \rel@kern{0.8}
    \overline{\rel@kern{-0.8}\macc@nucleus\rel@kern{0.2}}
    \rel@kern{-0.2}
  }
  \macc@depth\@ne
  \let\math@bgroup\@empty \let\math@egroup\macc@set@skewchar
  \mathsurround\z@ \frozen@everymath{\mathgroup\macc@group\relax}
  \macc@set@skewchar\relax
  \let\mathaccentV\macc@nested@a
  \macc@nested@a\relax111{#1}
  \endgroup
}
\makeatother

\let\oldcitet=\citet
\let\oldcitep=\citep
\renewcommand{\cite}[1]{\textcolor[rgb]{0,0,1}{\oldcitet{#1}}}
\renewcommand{\citet}[1]{\textcolor[rgb]{0,0,1}{\oldcitet{#1}}}
\renewcommand{\citep}[1]{\textcolor[rgb]{0,0,1}{\oldcitep{#1}}}

\renewcommand{\leq}{\leqslant}

\renewcommand{\geq}{\geqslant}

\begin{document}

\maketitle

\baselineskip0.6cm

\vspace{-0.6cm}

\begin{abstract}
We prove the existence of viscosity solutions for fractional semilinear elliptic PDEs on open balls with bounded exterior condition in dimension $d\geq 1$. Our approach relies on a tree-based probabilistic representation based on a $(2s)$-stable branching processes for all $s\in (0,1)$, and our existence results hold for sufficiently small exterior conditions and nonlinearity coefficients. In comparison with existing approaches, we consider a wide class of polynomial nonlinearities without imposing upper bounds on their maximal degree or number of terms. Numerical illustrations are provided in large dimensions. 
\end{abstract}


        
\baselineskip0.7cm

\parskip-0.1cm

\section{Introduction}
     Fully nonlinear Dirichlet problems
     for nonlocal operators have been studied in 
     \cite{bony} 
     by semi-group methods
     and in \cite{barles2}
      by the Perron method,
      for particular types of nonlinearities.
 For $d\geq 1$, let
      $$
 \Delta_s u = - ( - \Delta )^s u
 = \frac{4^s \Gamma (s + d/2 )}{\pi^{d/2}|\Gamma (- s)|}
 \lim_{r\rightarrow 0^+}
 \int_{\mathbb{R}^d \setminus B(x,r)} \frac{u(  \ \! \cdot +z)-u(z)}{|z|^{d+2s}}dz,
 $$
 denote the fractional Laplacian
 with parameter $s\in (0,1)$, 
 see, e.g., \cite{tendef},
      where $\Gamma (p) : = \displaystyle
   \int_0^\infty e^{-\lambda x} \lambda^{p-1} d\lambda$ 
   is the gamma function
   and $B(x,r)$ denotes the open ball of radius $r>0$ centered at $x\in \real^d$. 
   Consider the following nonlinear elliptic PDE 
   on an open set $\mathcal{O}$ in $\real^d$,    with 
      fractional Laplacian
      of the form  
 \begin{equation}
   \label{eq:1}
\begin{cases}
  \displaystyle
  \Delta_s u(x) +
 f(x,u(x))
  = u(x), \qquad x=(x_1,\ldots , x_d) \in \mathcal{O}, 
  \medskip
  \\
u(x) = \phi(x), \qquad x \in \mathbb{R}^d \setminus \mathcal{O},
\end{cases}
\end{equation}
 with polynomial non-linearities of the form
 $$
 f(x,y) = \sum\limits_{l \in \mathcal{L}} c_l(x) y^{l},
 \qquad
 x\in \mathcal{O},
 \quad
 y\in \real, 
 $$
 where 
\begin{itemize}
\item $\mathcal{L}$ is a subset of
 $\mathbb{N}$,
   \item
 $c_l(x)$, $l \in {\cal L}$,
 are bounded continuous functions on $\mathcal{O}$, and 
\item $\phi : \mathbb{R}^d \rightarrow \mathbb{R}$ is a Lipschitz function
 bounded on $\real^d \setminus \mathcal{O}$. 
 \end{itemize}    
      Two types of solutions,
      namely weak solutions, see Definition~2.1 in \cite{ros-oton},
      and viscosity solutions, see \cite{servadei2014},  
      are generally considered
      for elliptic PDEs. 
      While they coincide in many situations, see Remark~2.11 in
      \cite{ros-oton}, they involve different tools. 
      
      \medskip

Weak solutions can be obtained by the Riesz representation theorem
      or on the Lax-Milgram theorem, after rewriting the problem in its variational
      formulation, see \cite{ros-oton2016}, \cite{felsinger}. 
 The study of existence of viscosity solutions by the Perron method,
 see \cite{barles2}, 
 does not allow for general
 polynomial 
 non-linearities as in \eqref{eq:1},
 see conditions~(A1)-(A6) therein. 
 This restriction can be overcome by 
 semi-group methods \cite{bony}
 or using branching diffusion processes, see \cite{inw},
 or superprocesses, see \cite{LGBroSna}. 

      \medskip

      Under strong conditions on the nonlinearity $f(x,y)$, namely
      $|f(x,y)|\leq a_1 + a_2 |y|^{q-1}$ for some $q\in (2,2n/(n-2s))$
      and
      $\lim_{y\to \infty} f(x,y)/y=0$,
      existence of nontrivial weak solutions
     for problems of the form $ \Delta_s u(x) + f(x,u) =0$
     on an open bounded set $\mathcal{O}$ with $u=0$ on $\real^d \setminus \mathcal{O}$
     has been obtained in \cite{servadei} using the mountain pass theorem.
   In \cite{servadei2014}, existence of viscosity solutions
   has been proved for problems of the form  
   $ \Delta_s u(x) + f(x) =0$ with
   $u=\phi$ on $\real^d \setminus \mathcal{O}$ 
   under smoothness assumptions on 
   $f,\phi$,
   see also \cite{felsinger} and \cite{mou} for the existence of
   viscosity solutions, resp. weak solutions, with nonlocal operators. 

\medskip

On the other hand, a large part of the literature on fractional PDEs
with nonlinearities is devoted to proving the nonexistence of trivial solutions
when the initial datum $\phi$ vanishes
 outside $B(0,R)$ and $c_0(x) = 0$ in \eqref{eq:1}, see, e.g.,
 \cite{alves}, \cite{steglinski}, \cite{correia}. 
 This also includes the method of moving spheres, see e.g. \cite{pablo}, 
\cite{fall}, and the use of the Pohoazev identity for the fractional Laplacian,
 see \cite{ro2012pohone,ro2014poho}.
 We note that in this setting, 
 our existence results only yield the null function as a solution,
 which however does not contradict the nonexistence of nontrivial solutions.

   \medskip
   
   In this paper, we consider existence of 
   (continuous) viscosity solutions
   for fractional elliptic problems of the form \eqref{eq:1}
   using a large class of semilinearities $f(x,y)$
   without imposing upper bounds on their maximal
   degree or number of terms.
   We denote by
   $\widebar{B}(0,R)$ the
 closed  ball of radius $R>0$ in $\mathbb{R}^d$, 
 and for $\mathcal{O}$ an open subset of $\real^d$ we
 will use the fractional Sobolev space 
  $$H^s( \mathcal{O} ):=
  \left\{ u \in L^2( \mathcal{O} ) \ : \
  \frac{|u(x)-u(y)|}{|x-y|^{s+d/2}} \in L^2 (\mathcal{O}^2 )\right\}, 
$$
  with $d\geq 1$ and $s\in (0,1)$.
  Our main results can be stated as follows, with $\mathcal{O}=B(0,R)$. 
\begin{theorem}
\label{t2}
 Assume that $|\phi|_{L^\infty (B^c(0,R))} < \infty$ and
 $\displaystyle
 \sum_{l\in \mathcal{L}} |c_l|_{L^\infty (B(0,R))} < \infty$
 are both sufficiently small.
 Then, the nonlinear PDE~\eqref{eq:1}
 admits a (continuous) viscosity solution on $\mathcal{O}=B(0,R)$. 
\end{theorem}
We note from the proof of Theorem~\ref{t2} that
it suffices in particular to have
$ | \phi|_{L^\infty (B^c(0,R))} \leq 1$ and $\sum_{l \in \mathcal{L}} |c_l|_{L^\infty (B(0,R))} \leq 1$
for its conclusion to hold. 
\begin{theorem}
\label{t3}
 Assume that ${\cal L}$ is finite and that 
 $|\phi|_{L^\infty (B^c(0,R))} < \infty$ and
 $|c_l|_{L^\infty (B(0,R))} < \infty$,
 $l\in \mathcal{L}$. 
 Then, the nonlinear PDE~\eqref{eq:1}
 admits a (continuous) viscosity solution on $\mathcal{O}=B(0,R)$
 for sufficiently small $R>0$. 
\end{theorem}
We note in Proposition~\ref{p3} that
 if $\phi$ belongs to $H^{2s} (\real^d)$,  the viscosity solution 
 obtained in Theorems~\ref{t2}-\ref{t3}
 is also a weak solution,
 and that it is the only weak and viscosity solution
 of \eqref{eq:1}. 

    \medskip 

    Theorems~\ref{t2}-\ref{t3} will be proved through a probabilistic representation of PDE solutions using branching stochastic processes.
    Stochastic branching processes 
for the representation of PDE solutions 
have been introduced by \cite{skorohodbranching},
\cite{inw}, and have been used to prove blow-up and existence of
 solutions for parabolic PDEs
 in \cite{N-S}, \cite{lm}, \cite{penent}. 

\medskip

\indent
 This branching argument has been recently applied in
 \cite{labordere} to the treatment of
 parabolic PDEs with polynomial gradient nonlinearities, 
 see \cite{claisse} for the elliptic case.
 In this approach, gradient terms
 are associated to tree branches to which
 a Malliavin integrations by parts is applied.
 In \cite{penent}, 
 this approach has been extended to
 semilinear parabolic PDEs with pseudo-differential
 operators of the form $-\eta(-\Delta /2)$ 
 and fractional Laplacians,
 using a random tree $\mathcal{T}_x $ starting at $x\in \real^d$
 and carrying a symmetric $(2s)$-stable process.
 In the absence of gradient nonlinearities,
 the tree-based approach has been recently implemented
 for nonlocal semilinear parabolic PDEs in \cite{belak}.

\medskip 
 
 In what follows, we will apply
 this probabilistic representation approach to the setting
 of elliptic PDEs with fractional Laplacians. 
 PDE solutions
 will be constructed as the expectation
 $u(x) = \mathbb{E}[\mathcal{H}(\mathcal{T}_{x})]$ 
 of a
 random functional $ \mathcal{H}(\mathcal{T}_x)$ of the underlying branching process,
 $x\in \real^d$, 
 which yields a probabilistic representations for
the solutions of a wide class of semilinear elliptic PDEs of the form
\eqref{eq:1}.
Sufficient conditions for the representation of classical 
solutions of \eqref{eq:1} as 
 $u(x) = \mathbb{E}[\mathcal{H}(\mathcal{T}_{x})]$ 
 are obtained in Proposition~\ref{p2}. 

\medskip

As we are dealing with continuous viscosity solutions, we need to ensure
that the random variable $\mathcal{H}(\mathcal{T}_x)$ is sufficiently integrable, 
so that the expected value
$\mathbb{E}[\mathcal{H}(\mathcal{T}_{x})]$ is a continuous function of $x\in \real^d$,
see Lemma~\ref{pl1}.  
For this, in Lemma~\ref{l1} we show, using results of \cite{kyprianou2018},  
that the exit time of a stable process from the ball
$B(0,R)$ is almost surely continuous with respect to
its initial condition $x\in B(0,R)$.
This allows us to
show the uniform integrability
 of $(\mathcal{H}(\mathcal{T}_{x}))_{x \in B(0,R)} $ 
using 
the fractional Laplacian $ \Delta_s =-(-\Delta)^s $ and its
associated stable process.
We note that the result of Lemma~\ref{pl1} 
 may be extended to non spherical
 domains for which Lemma~\ref{l1} is satisfied. 
 
\medskip 

 In Section~\ref{s5}
 provide a numerical implementation of our existence results
 using Monte Carlo simulations
 for nonlinear fractional PDEs
 in dimension up to 100.
 We note that the tree-based Monte Carlo method   
 allows us to solve large dimensional problems, 
 whereas the application of deterministic finite difference
 methods is generally restricted to one dimension and 
 their extension to higher dimensions still remains a challenge,
 see e.g. \cite{huang-oberman} in the linear case.  
 
\medskip

\indent This paper is organized as follows.
 The description of the branching mechanism is
 presented in Section~\ref{s2}.
 In Section~\ref{s3} we state and prove 
 Theorem~\ref{t1.0} 
 which gives the probabilistic representation of the solution and
 its partial derivatives.
 Finally, in Section~\ref{s5} 
 we present numerical simulations to illustrate the method on specific examples.
 
 \medskip
  
 The C~codes used to plot 
 Figures~\ref{fig1} to \ref{fig3} 
 are available at
 
 \centerline{\url{https://github.com/nprivaul/semilinear_fractional_elliptic}.}
 
\section{Probabilistic representation of elliptic PDE solutions}
\label{s2}
This section describes the
probabilistic representation for the solution of \eqref{eq:1}, using a branching mechanism
 giving the solution of \eqref{eq:1}
 as the expectation of a multiplicative functional defined
on a random tree structure. 
 The probabilistic representations of 
 Theorems~\ref{t2}-\ref{t3} use a functional
 on a random branching process driven
 by a stable L\'evy process $(X_t)_{t\in \real_+}$ 
 such that
 $\mathbb{E}\big[e^{ i \xi X_t}\big] = e^{ - t | \xi |^{2s }}$,
 where $|\xi |$ denotes the Euclidean norm of 
 $\xi\in \real^d$, $t \geq 0$,
 so that
 the infinitesimal generator 
 of $(X_t)_{t\in \real_+}$
 is
 the fractional Laplacian
 $\Delta_s = - ( - \Delta )^s$,
 $s \in (0,1]$. 
\subsubsection*{Random tree} 
 Given $\rho: \mathbb{R}^+ \rightarrow (0,\infty )$ a
 probability density function on $\real_+$,
 consider a probability mass function
$(q_{l})_{l\in {\cal L}}$ on ${\cal L}$
 with $q_{l} > 0$,
$l\in {\cal L}$, and
$\sum_{ l \in {\cal L}} l q_{l} < \infty$.
 In addition, we consider, on a probability space $(\Omega,\mathcal{F},\mathbb{P})$, 
\begin{itemize}
\item an {\it i.i.d.} family $(\tau^{i,j})_{i,j\geq 1}$ of random variables
 with distribution $\rho (t)dt$ on $\real_+$,
\item an {\it i.i.d.} family $(I^{i,j})_{i,j\geq 1}$ of discrete
  random variables, with
  $$
  \mathbb{P} \big( I^{i,j}=l \big) = q_l >0,
  \qquad l \in {\cal L},
  $$
\item 
  an independent family $\big(X^{i,j}\big)_{i,j\geq 1}$
  of symmetric $(2s)$-stable processes, 
\end{itemize}
 where the sequences $(\tau^{i,j})_{i,j\geq 1}$, $(I^{i,j})_{i,j\geq 1}$ and
 $\big(X^{i,j}\big)_{i,j\geq 1}$ are assumed to be mutually independent.

 \subsubsection*{Branching process} 

We consider a branching process starting
from a particle $x\in B(0,R)$ with label $\widebar{1}=(1)$, 
which evolves according to
the process $X_{s,x}^{\widebar{1}} = x + X_{s}^{1,1}$,
$s \in [0,T_{\widebar{1}} ]$
with
$T_{\widebar{1}} := \min \big( \tau^{1,1} , \tau^B (x) \big)$, where
$$
\tau^B (x) := \inf\left\{ t \geq 0 \ : \ x + X_t^{1,1} \not\in B(0,R) \right\}
$$ 
 denotes the first exit time of
 $\big( x + X_s^{1,1}\big)_{s\in \real_+}$ 
from $B(0,R)$ 
 after starting from $x\in B(0,R)$. 
 Note that
 by (1.4) in \cite{bogdanbarrier} we have 
 $\E \big[ \tau^B (x) \big] < \infty$, and therefore
 $\tau^B(x)$ is almost surely finite for all $x\in B(0,R)$.
 On the other hand, although $\tau^B (x)$ depends on $(1,1)$,
 for the sake of clarity we will omit this information 
 in the sequel. 
 
\medskip 

If $\tau^{1,1}<\tau^B(x)$, the process branches at time $\tau^{1,1}$
into new independent copies of
$(X_t)_{t \in \real_+}$, each of them started at
$X_{x,\tau^{1,1}}^{\widebar{1}}$.
Based on the values of $I^{1,1} = l \in {\cal L}$,
a family of $l$ of new branches is created with the probability $q_l$,
 where
\begin{itemize}
\item the first $l_0$ branches are indexed by $(1,1),(1,2),\ldots ,(1,l_0)$,
\item the next $l_1$ branches
  are indexed by $(1,l_0+1),\ldots ,(1,l_0 + l_1)$, and so on.
\end{itemize}
Each new particle then follows independently
the same mechanism as the first one, and
every branch stops when it leaves the domain $B(0,R)$.
Particles at generation $n\geq 1$ are assigned a label of the form
 $\widebar{k} = (1,k_2,\ldots ,k_n) \in \mathbb{N}^n$,
and their parent is labeled $\widebar{k}- := (1,k_2,\ldots ,k_{n-1})$.
The particle labeled $\widebar{k}$ is born at time $T_{\widebar{k}-}$
and its lifetime $\tau^{n,\pi_n(\widebar{k})}$ is the element of index
$\pi_n\big(\widebar{k}\big)$ in the {\it i.i.d.} sequence
$(\tau^{n,j})_{j\geq 1}$ with
 $\tau^{1,\pi_1(\widebar{1})} = \tau^{1,1}$, 
defining an injection
$$
\pi_n:\mathbb{N}^n \to \mathbb{N}, \qquad n\geq 1, 
$$ 
such that $\pi_1(1)=1$.
The random evolution of particle $\widebar{k}$
is given by
\begin{equation}
\nonumber 
X_{s,x}^{\widebar{k}} := X^{\widebar{k}-}_{T_{\widebar{k}-},x}+X_{s-T_{\widebar{k}-}}^{n,\pi_n(\widebar{k})},
\qquad s\in [T_{\widebar{k}-},T_{\widebar{k}}],
\end{equation} 
where $T_{\widebar{k}} := T_{\widebar{k}-} +
\min \big(
\tau^{n,\pi_n(\widebar{k})} , \tau^B \big( {X^{\widebar{k}-}_{T_{\widebar{k}-},x}}\big) \big)$,
$\widebar{k} \in \mathbb{N}^n$, $n\geq 2$, and
$$
\tau^B \big( {X^{\widebar{k}-}_{T_{\widebar{k}-},x}}\big):= \inf\Big\{ t \geq 0 \ : \ X^{\widebar{k}-}_{T_{\widebar{k}-},x}+X_t^{n,\pi_n(\widebar{k})} \not\in B(0,R) \Big\}. 
$$
 Given
$\widebar{k} = (1,k_2,\ldots ,k_n) \in \mathbb{N}^n$,
if $\tau^{n,\pi_n(\widebar{k})} < \tau^B\big( {X^{\widebar{k}-}_{T_{\widebar{k}-},x}}\big)$,
 we draw a sample
 $I_{\widebar{k}} : = I^{n,\pi_n(\widebar{k})} = l$
 of $I^{n,\pi_n(\widebar{k})}$,
 and the particle $\widebar{k}$ branches into
 $\big|I^{n,\pi_n(\widebar{k})}\big|=l$ offsprings at generation
  $n+1$,
 which are indexed by $(1,\ldots ,k_n,i)$, $i=1,\ldots ,\big|I^{n,\pi_n(\widebar{k})}\big|$.

 \medskip 

 The set of particles dying inside $B(0,R)$ is denoted by $\mathcal{K}^{\circ}$,
 whereas those dying outside form a set denoted by $\mathcal{K}^{\partial}$. The particles of $n$-th generation, $n\geq 1$,
 will be denoted by $\mathcal{K}_n^\circ$
 (resp. $\mathcal{K}_n^\partial$) if they die inside the domain (resp. outside). We also
  define the filtration $(\mathcal{F}_n)_{n\geq 1}$ as 
$$
\mathcal{F}_n := \sigma\left(T_{\widebar{k}},I_{\widebar{k}},X^{\widebar{k}},\widebar{k} \in \bigcup_{i=1}^n \mathbb{N}^i\right), \qquad n \geq 1. 
$$

\begin{definition}
 When started from a position $x\in \real^d$,  
 the above construction yields 
 a branching process called a random tree rooted at $x$, and 
 denoted by $\mathcal{T}_x$. 
\end{definition}
 The tree $\mathcal{T}_{x}$
 will be used for the stochastic representation of the solution
 $u(x)$ of the PDE~\eqref{eq:1}.
  The next table summarizes the notation introduced so far.

 \medskip

 \bigskip

\begin{center}
 \begin{tabular}{||l | c||}
 \hline
 Object & Notation \\ [0.5ex]
 \hline\hline
 Initial position & $x$ \\
 \hline
 Tree rooted at $x$ &  $\mathcal{T}_x$ \\  \hline
 Particle (or label) of generation $n\geq 1$ & $\stackrel{}{\widebar{k}}=(1,k_2,\ldots ,k_n)$\\
 \hline
 First branching time & $T_{\widebar{1}}$\\
 \hline
 Lifespan of a particle & $T_{\widebar{k}} - T_{\widebar{k}-}$ \\
 \hline
Birth time of a particle $\widebar{k}$ & $T_{\widebar{k}-}$ \\
\hline
Death time of a particle $\widebar{k} \in \mathcal{K}^\circ$& $T_{\widebar{k}} = \T_{\widebar{k}-} + \tau^{n,\pi_n(\widebar{k})}$ \\
\hline
Death time of a particle $\widebar{k} \in \mathcal{K}^\partial$& $T_{\widebar{k}} = T_{\widebar{k}-} + \tau^B\big( {X^{\widebar{k}-}_{T_{\widebar{k}-},x}}\big)$ \\
\hline
Position at birth & $X^{\widebar{k}}_{T_{\widebar{k}-},x}$\\
\hline
Position at death & $X^{\widebar{k}}_{T_{\widebar{k}},x}$ \\
\hline
Exit time starting from $x$ & $\tau^B (x):= \inf\left\{ t \geq 0 \ : \ x + X_t \not\in B(0,R) \right\}$
\\ 
\hline
\end{tabular}
\end{center}

\bigskip
\noindent
To represent the structure
of the tree we use the following conventions, in which
different colors represent different ways of branching:

\medskip 
  
\tikzstyle{level 1}=[level distance=4cm, sibling distance=4cm]
\tikzstyle{level 2}=[level distance=5cm, sibling distance=3cm]

\begin{center}
\resizebox{0.45\textwidth}{!}{
\begin{tikzpicture}[scale=0.9,grow=right, sloped][H]
\node[ellipse split,draw,purple,text=black,thick]{Time \nodepart{lower} Position}
    child {
        node[ellipse split,draw,purple,text=black,thick]{Time \nodepart{lower} Position}
        child{
        node[ellipse,draw,thick]{...}
        edge from parent
        node[above]{Label}
        }
        child{
        node[ellipse,draw,thick]{...}
        edge from parent
        node[above]{Label}
        } 
        edge from parent
        node[above]{Label}
    }
    child {
        node[ellipse split,draw,cyan,text=black,thick]{Time \nodepart{lower} Position}
        edge from parent
        node[above]{Label}
    };
\end{tikzpicture}
}
\end{center}

\vskip-0.2cm

\noindent
 Specifically, let us draw a tree sample for the PDE 
$$
 \Delta_s u (x) + c_0 (x) + c_{0,1}(x) u^2 (x) = 0
$$
in dimension $d=1$. For this tree, there are two types of branching: we can either branch into no branch at all (which is represented in cyan color), or into two branches. 
The black color is used for leaves, namely the particles that leave the domain $B(0,R)$. 
\\
 
\begin{center}
\resizebox{0.90\textwidth}{!}{
\begin{tikzpicture}[scale=1.1,grow=right, sloped, yscale = 0.85]
\node[ellipse split,draw,blue,thick]{$0$ \nodepart{lower} $x$}
    child {
        node[ellipse split,draw,purple,text=black,thick] {$T_{\widebar{1}}$ \nodepart{lower} $X^{\widebar{1}}_{T_{\widebar{1}},x} $}        
            child {
                node[ellipse split,draw,purple,text=black,thick] {$T_{(1,2)}$ \nodepart{lower} $X^{(1,2)}_{T_{(1,2)},x} $} 
                child{
                node[ellipse split,draw,cyan,text=black,thick, right=1cm]{$T_{(1,2,2)}$ \nodepart{lower} $X^{(1,2,2)}_{T_{(1,2,2)},x} $}
                edge from parent
                node[above]{~~~~~$(1,2,2)$~~~~~}
                }
                child{
                node[ellipse split,draw,thick, right=0.1cm]{$T_{(1,2,1)}:=T_{(1,2)}+\tau^B\big( {X^{(1,2,2)}_{T_{(1,2,2)},x}}\big)$ \nodepart{lower} $X^{(1,2,1)}_{T_{(1,2,1)},x} $}
                edge from parent
                node[above]{~~~~~$(1,2,1)$~~~~~}
                }
                edge from parent
                node[above] {$(1,2)$}
            }
            child {
                node[ellipse split,draw,cyan,text=black,thick] {$T_{\widebar{1}} + T_{(1,1)}$ \nodepart{lower} $X^{(1,1)}_{T_{(1,1)},x}  $}
                edge from parent
                node[above] {$(1,1)$}
            }
            edge from parent 
            node[above] {$\widebar{1}$}
    };
\end{tikzpicture}
}
\end{center}
\noindent 
 In the above example we have
 $\mathcal{K}^{\circ}= \big\{\widebar{1}, (1,1) , (1,2) ,(1,2,2)\big\}$
 and $\mathcal{K}^{\partial} = \{(1,2,1)\}$.

\subsubsection*{Representation of PDE solutions}
Given $x\in \mathbb{R}^d$, 
we consider the functional $\mathcal{H}$ 
  of the random tree $\mathcal{T}_x$ defined as 

\begin{equation}
\label{djsdiuoa} 
        {\mathcal{H} (\mathcal{T}_x) :=
          \prod_{\widebar{k} \in \mathcal{K}^{\circ}} \frac{e^{-{\Delta T_{\widebar{k}}}}c_{I_{\widebar{k}}}\big(X^{\widebar{k}}_{T_{\widebar{k}},x}\big)}{q_{I_{\widebar{k}}}\rho(
\Delta T_{\widebar{k}} 
)} \prod_{\widebar{k} \in \mathcal{K}^{\partial}} \frac{ e^{-{\Delta T_{\widebar{k}}}} \phi\big(X^{\widebar{k}}_{T_{\widebar{k}},x}\big) }{\widebar{F}(
            \Delta T_{\widebar{k}}
             )}},
\end{equation} 
where
$\Delta T_{\widebar{k}} : = T_{\widebar{k}} - T_{\widebar{k}-}$,
$\widebar{k} \in \mathcal{K}$,
and 
$\widebar{F}(t) : = 1 - \mathbb{P}(T_{\widebar{1}} \leq t )$.
The next proposition provides a probabilistic representation for the
solution of \eqref{eq:1} as the expected value of
 the functional $\mathcal{H} (\mathcal{T}_x)$. 
\begin{prop}
  \label{p2}
  Assume that the PDE~\eqref{eq:1} admits a classical solution $u \in H^{2s} (B(0,R)) \cap \mathcal{C}\big(\widebar{B}(0,R)\big)$, 
    such that the sequence 
$$
         {\mathcal{H}_n(\mathcal{T}_x) :=
          \prod_{\widebar{k} \in \bigcup_{i=1}^n \mathcal{K}_i^{\circ}} \frac{e^{-{\Delta T_{\widebar{k}}}}c_{I_{\widebar{k}}}\big(X^{\widebar{k}}_{T_{\widebar{k}},x}\big)}{q_{I_{\widebar{k}}}\rho(
            \Delta T_{\widebar{k}}
            )}  \prod_{\widebar{k} \in \bigcup_{i=1}^n  \mathcal{K}_i^{\partial}} \frac{ e^{-{\Delta T_{\widebar{k}}}} \phi\big(X^{\widebar{k}}_{T_{\widebar{k}},x}\big) }{\widebar{F}(
            \Delta T_{\widebar{k}}
            )}} \prod_{\widebar{k} \in \mathcal{K}_{n+1}} u\big(X_{T_{\widebar{k}-}^{\widebar{k}} }\big)
$$ 
        is uniformly integrable
        in $n \geq 1$.  
Then we have $u(x) = \mathbb{E}[\mathcal{H}(\mathcal{T}_{x})]$,
$x\in B(0,R)$. 
\end{prop} 
\begin{Proof}
  Applying the It\^o-Dynkin formula to the process
  $\big(e^{-t}u\big(X_{t,x}^{\widebar{1}}\big)\big)_{t\in \real_+}$
   on the time interval $\big[0,\tau^B (x)\big]$, we find 
 $$
 \mathbb{E}\big[e^{-\tau^B(x)}u\big(X_{\tau^B(x),x}^{\widebar{1}} \big)\big] = u(x)
 + \mathbb{E}\bigg[  \int_0^{\tau^B(x)} e^{-t}
   \Delta_s
    u\big(X_{t,x}^{\widebar{1}} \big) dt - \int_0^{\tau^B(x)} e^{-t} u\big(X_{t,x}^{\widebar{1}} \big) dt\bigg], 
$$
 which implies that $u(x)$ can be represented as
 $$
 u(x) = \mathbb{E}\bigg[ e^{-\tau^B(x)} u\big(X_{\tau^B(x),x}^{\widebar{1}} \big) + \int_0^{\tau^B(x)} e^{-t}f\big(X_{t,x}^{\widebar{1}} ,u\big( X_{t,x}^{\widebar{1}} \big)\big) dt \bigg], 
$$
 since $u$ solves \eqref{eq:1}.  
Therefore, since $T_{\widebar{1}}$ has the
probability density $\rho$ and 
is independent of $\big(X_{s,x}^{\widebar{1}}\big)_{s\in \real_+}$, 
 from the boundary condition $u(x)=\phi(x)$, $x\in \mathbb{R}^d \setminus B(0,R)$,
 we have 
\begin{eqnarray}
  \nonumber
  u(x) & = &  
  \mathbb{E}\left[
    \mathbb{E}\left[
      e^{-\tau^B(x)} u\big(X_{\tau^B(x),x}^{\widebar{1}}\big) + \int_0^{\tau^B(x)} e^{-t}
      f\big(X_{t,x}^{\widebar{1}},u\big(X_{t,x}^{\widebar{1}}\big)\big)
       dt
     \ \! \bigg|
      \ \!
      \big(X_{t,x}^{\widebar{1}}\big)_{t\in \real_+}
      \right]\right]
\\
\nonumber
    &= & 
  \mathbb{E}\left[
      \mathbb{E}\left[
        \frac{e^{-\tau^B(x)}}{\widebar{F}(\tau^B(x))}
         \phi\big(X_{\tau^B(x),x}^{\widebar{1}}\big) \mathbbm{1}_{\{ T_{\widebar{1}}= \tau^B(x) \} }
 \ \! \bigg|
      \ \!
      \big(X_{t,x}^{\widebar{1}}\big)_{t\in \real_+}
      \right]
        \right]
\\
  \nonumber
  &  &
  + \mathbb{E}\left[
      \mathbb{E}\left[
        \int_0^{\tau^B(x)} \frac{e^{-t}}{\rho (t) }
        f\big(X_{t,x}^{\widebar{1}},u\big(X_{t,x}^{\widebar{1}}\big) \big) 
          \rho (t) dt
     \ \! \bigg|
      \ \!
      \big(X_{t,x}^{\widebar{1}}\big)_{t\in \real_+}
      \right]\right]
 \\
\nonumber
    & = &  
\mathbb{E}\left[ \frac{e^{-\tau^B(x)}}{\widebar{F}(\tau^B(x))} \phi\big(
  X_{\tau^B(x),x}^{\widebar{1}}\big) \mathbbm{1}_{\{ T_{\widebar{1}}= \tau^B (x)\} }
  +
  \frac{e^{-T_{\widebar{1}}}}{\rho(T_{\widebar{1}}) }
   f\big(X_{t,x}^{\widebar{1}},u\big(X_{t,x}^{\widebar{1}}\big)\big)
   \mathbbm{1}_{\{ T_{\widebar{1}} < \tau^B(x) \} }   \right]
 \\
\label{fkjdslf3} 
&= &  \mathbb{E}\left[  \frac{e^{- T_{\widebar{1}}}}{\widebar{F}( T_{\widebar{1}})} \phi\big(
  X_{\tau^B(x),x}^{\widebar{1}}
  \big) \mathbbm{1}_{\{ T_{\widebar{1}}= \tau^B(x) \} } +  \frac{e^{-T_{\widebar{1}}}}{\rho(T_{\widebar{1}}) } \frac{c_{I_{\widebar{1}}}\big(
    X_{T_{\widebar{1}},x}^{\widebar{1}}
    \big)}{q_{I_{\widebar{1}}}}
   u^{I_{\widebar{1}}} \big(X_{T_{\widebar{1}},x}^{\widebar{1}}\big)
    \mathbbm{1}_{ \{ T_{\widebar{1}} < \tau^B(x) \} }   \right], \qquad 
\end{eqnarray}
showing that $u(x) = \mathbb{E}[\mathcal{H}_1 (\mathcal{T}_x ) ]$,
$x\in B(0,R)$, since $\mathcal{K}_1 = \{\widebar{1} \}$,  
$$
        {\mathcal{H}_1(\mathcal{T}_x) =
          \prod_{\widebar{k} \in \mathcal{K}_1^{\circ}} \frac{e^{-{\Delta T_{\widebar{k}}}}c_{I_{\widebar{k}}}\big(X^{\widebar{k}}_{T_{\widebar{k}},x}\big)}{q_{I_{\widebar{k}}}\rho(
            \Delta T_{\widebar{k}}
            )}  \prod_{\widebar{k} \in \mathcal{K}_1^{\partial}} \frac{ e^{-{\Delta T_{\widebar{k}}}} \phi\big(X^{\widebar{k}}_{T_{\widebar{k}},x}\big) }{\widebar{F}(
            \Delta T_{\widebar{k}}
            )}} \prod_{\widebar{k} \in \mathcal{K}_2} u\big(X_{T_{\widebar{k}-}^{\widebar{k}} }\big),
$$ 
 and $X_{T_{\widebar{1}},x}^{\widebar{1}} = X_{T_{\widebar{k}-},x}^{\widebar{k}}$. 
        Repeating the above argument after starting from
$X_{T_{\widebar{1}},x}^{\widebar{1}} = X_{T_{\widebar{k}-},x}^{\widebar{k}}$
 instead of $x$,
 $\widebar{k}\in \mathcal{K}_2$,
 and using the independence of
$\big(X_{T_{\widebar{k}},x}^{\widebar{k}}\big)_{\widebar{k}\in \mathcal{K}_2}$ 
 given $\mathcal{F}_1$, we find 
\begin{align*} 
 & 
    u^{I_{\widebar{1}}}\big(X_{T_{\widebar{1}},x}^{\widebar{1}}\big) = 
  \prod_{\widebar{k}\in \mathcal{K}_2} 
  u\big(X_{T_{\widebar{k}-},x}^{\widebar{k}}\big)
  \\
   & =  
  \prod_{\widebar{k}\in \mathcal{K}_2}
  \mathbb{E}\left[ \frac{e^{- T_{\widebar{k}}}}{\widebar{F}( T_{\widebar{k}})} \phi\big(
   X_{T_{\widebar{k}},x}^{\widebar{k}}
   \big) \mathbbm{1}_{\{
     X_{T_{\widebar{k}},x}^{\widebar{k}}
     \notin B(0,R)\}} + \frac{e^{-T_{\widebar{k}}}}{\rho(T_{\widebar{k}}) } \frac{c_{I_{\widebar{k}}}\big(
     X_{T_{\widebar{k}},x}^{\widebar{k}}
\big)}{q_{I_{\widebar{k}}}}
   u^{I_{\widebar{k}}} \big(X_{T_{\widebar{k}},x}^{\widebar{k}}\big)
   \mathbbm{1}_{\{ X_{T_{\widebar{k}},x}^{\widebar{k}} \in B(0,R) \} }
   \ \! \bigg| \ \! \mathcal{F}_1 \right]
   \\
   & =  
   \mathbb{E}\left[
     \prod_{\widebar{k}\in \mathcal{K}_2}
     \left(
     \frac{e^{- T_{\widebar{k}}}}{\widebar{F}( T_{\widebar{k}})} \phi\big(
   X_{T_{\widebar{k}},x}^{\widebar{k}}
   \big) \mathbbm{1}_{\{
     X_{T_{\widebar{k}},x}^{\widebar{k}}
     \notin B(0,R)\}} + \frac{e^{-T_{\widebar{k}}}}{\rho(T_{\widebar{k}}) } \frac{c_{I_{\widebar{k}}}\big(
     X_{T_{\widebar{k}},x}^{\widebar{k}}
\big)}{q_{I_{\widebar{k}}}}
   u^{I_{\widebar{k}}} \big(X_{T_{\widebar{k}},x}^{\widebar{k}}\big)
   \mathbbm{1}_{\{ X_{T_{\widebar{k}},x}^{\widebar{k}} \in B(0,R) \} }
   \right)
   \ \! \bigg| \ \! \mathcal{F}_1 \right]. 
\end{align*} 
Plugging this expression in \eqref{fkjdslf3} above and using the tower property
of the conditional expectation 
yields $u(x)= \mathbb{E}[\mathcal{H}_2(\mathcal{T}_x )]$,
and repeating this process inductively leads to
$u(x)= \mathbb{E}[\mathcal{H}_n(\mathcal{T}_x )]$.
From the uniform integrability of 
$(\mathcal{H}_n(\mathcal{T}_x ))_{n\geq 1}$
        and the fact that the tree goes extinct almost surely,
        we conclude to $u(x) = \mathbb{E}[\mathcal{H}(\mathcal{T}_{x})]$
        after letting $n$ tend to infinity. 
\end{Proof}
\section{Proofs of existence results} 
\label{s3}
In Proposition~\ref{p1} we start by showing 
that  existence of solutions holds 
when $u(x):= \mathbb{E}[\mathcal{H}(\mathcal{T}_{x})]$
is continuous in $ x\in \widebar{B}(0,R)$. 
 We then prove in Lemma~\ref{pl1} that continuity of
 $\mathbb{E}[\mathcal{H}(\mathcal{T}_{x})]$
 in $ x\in \widebar{B}(0,R)$ holds
 when $(\mathcal{H}(\mathcal{T}_{x}))_{x\in B(0,R)}$
is uniformly integrable. 
The proofs of Theorems~\ref{t2}-\ref{t3} 
are then completed by showing that
 $(\mathcal{H}(\mathcal{T}_{x}))_{x \in B(0,R)} $ is uniformly bounded
in $L^p (B(0,R))$ for some $p>1$, implying the required uniform integrability. 
\begin{prop}{}
\label{p1}
Suppose that
$\E\big[ \big| \mathcal{H}(\mathcal{T}_{x})\big| \big]<\infty$ for all
$x\in B(0,R)$,
and that $u(x):= \mathbb{E}[\mathcal{H}(\mathcal{T}_{x})]$ is
continuous on $\widebar{B}(0,R)$.
Then, $u$ is a viscosity solution of the PDE~\eqref{eq:1}.
\end{prop}
\begin{Proof}
  By conditioning with respect to
  $X^{\widebar{1}}_{T_{\widebar{1}},x}$ and $I_{\widebar{1}}$ 
  and using the fact that each
  offspring starts
  independent identically distributed
  stable branching processes, we have 
$$
\mathbb{E}\left[\prod_{i=0}^{I_{\widebar{1}}-1} \mathcal{H}\big(\mathcal{T}_{X^{\widebar{1}}_{T_{\widebar{1}}},x}\big)
  \ \bigg| \ X^{\widebar{1}}_{T_{\widebar{1}},x},I_{\widebar{1}} \right]
 \mathbbm{1}_{\{ T_{\widebar{1}} < \tau^B(x) \} } = u^{I_{\widebar{1}}}\big(X^{\widebar{1}}_{T_{\widebar{1}},x}\big)\mathbbm{1}_{\{ T_{\widebar{1}} < \tau^B(x) \}}, 
$$
 hence, since $T_{\widebar{1}}$ has density $\rho$, we get 
 \begin{eqnarray*}
   u(x) & = &
   \mathbb{E}[\mathcal{H}(\mathcal{T}_{x})]
   \\
    & = & 
   \mathbb{E}\left[ \frac{e^{- T_{\widebar{1}}}}{\widebar{F}( T_{\widebar{1}})} \phi\big(X^{\widebar{1}}_{\tau^B(x),x}\big) \mathbbm{1}_{\{ T_{\widebar{1}} = \tau^B(x) \} } + \frac{e^{-T_{\widebar{1}}}}{\rho(T_{\widebar{1}}) } \frac{c_{I_{\widebar{1}}}\big(X^{\widebar{1}}_{T_{\widebar{1}},x}\big) }{q_{{I_{\widebar{1}}}} } \prod_{i=0}^{I_{\widebar{1}}-1} \mathcal{H}\big(\mathcal{T}_{X^{\widebar{1}}_{T_{\widebar{1}}},x}\big) \mathbbm{1}_{\{ T_{\widebar{1}} < \tau^B(x) \} } \right]
   \\
    & = & 
   \mathbb{E}\left[ e^{-\tau_x^B} \phi\big(X^{\widebar{1}}_{\tau^B(x),x}\big) + \int_0^{\tau^B(x)} e^{-t} f\big(X^{\widebar{1}}_{t,x} ,u\big(X^{\widebar{1}}_{t,x}\big)\big) dt \right]
   \\
    & = & 
 \mathbb{E}\left[ e^{- \delta \wedge \tau_x^B} u\big(X^{\widebar{1}}_{\delta \wedge \tau^B(x),x}\big) + \int_0^{\delta \wedge \tau^B(x)} e^{-t} f\big(X^{\widebar{1}}_{t,x} ,u\big(X^{\widebar{1}}_{t,x}\big)\big) dt \right], 
\end{eqnarray*}
 for any $\delta > 0$, by the Markov property.
 It then follows from a classical argument that $u$ is a viscosity solution
 of the PDE~\eqref{eq:1}.
 Indeed, 
 let $ \xi \in \mathcal{C}^2(B(0,R))$ such that $x$ is a maximum point of $u-\xi$ and $u(x) = \xi(x)$. By the It\^o-Dynkin formula, we get 
$$
\mathbb{E}\big[e^{-\delta \wedge \tau^B(x)}\xi\big(X^{\widebar{1}}_{\delta \wedge \tau^B(x),x}\big)\big]
 = \xi(x)
 + \mathbb{E}\bigg[  \int_0^{\delta \wedge \tau^B(x)} e^{-t} \Delta_s \xi \big(X^{\widebar{1}}_{t,x}\big) dt - \int_0^{\delta \wedge \tau^B(x)} e^{-t} \xi \big(X^{\widebar{1}}_{t,x}\big) dt\bigg]. 
$$
Thus, since $u(x) = \xi(x)$ and $u\leq \xi$, we obtain 
$$
\mathbb{E}\left[ \int_0^{\delta \wedge \tau^B(x)} e^{-t}\big( 
  \Delta_s \xi \big(X^{\widebar{1}}_{t,x}\big)- \xi \big(X^{\widebar{1}}_{t,x}\big)+ f\big(X^{\widebar{1}}_{t,x} , u\big(X^{\widebar{1}}_{t,x}\big)\big)
  \big) dt \right] \geq 0. 
$$
Since $X_t^x$ converges in distribution to the constant $x\in \real^d$
as $t$ tends to zero,
it admits an  almost surely convergent subsequence,
hence by continuity and boundedness of $f( \ \! \cdot \ \! ,u( \ \! \cdot  \ \! ))$
together with the mean-value and dominated convergence theorems,
we have 
 $$
 \Delta_s \xi(x)  + f(x,\xi(x)) - \xi(x) \geq 0. 
$$
 We conclude that $u$ is a viscosity subsolution
 (and similarly a viscosity supersolution)
 of \eqref{eq:1}.
\end{Proof}
\begin{lemma}
\label{l1}
 The following statements hold true with probability one.
\begin{enumerate}[a)]
\item
 Let $x \in B(0,R)$.
 For $\P$-almost every $\omega \in \Omega $ there exists $r_0 (\omega )>0$ such that
$
\tau^B(y) = \tau^B(x), 
$
for all $y \in B(x,r_0 (\omega ))$.
\item
  Let $x \in \widebar{B}(0,R) \setminus B(0,R)$. 
  We have $\lim_{n\to \infty} \tau^B(x_n) = \tau^B(x)=0$ 
  almost surely for any sequence $(x_n)_{n\geq 0}$ in $B(0,R)$
  converging fast enough to $x$.
\end{enumerate}
 As a consequence, for any $x\in \widebar{B}(0,R)$ we have 
  \begin{eqnarray}
\label{eq:tau}
\mathbb{P}\Big(\lim_{n\rightarrow \infty} \tau^B (x_n) = \tau^B(x) \Big) = 1. 
\end{eqnarray}
 for any sequence $(x_n)_{n\in \inte}$ in $B(0,R)$
 converging fast enough to $x \in \widebar{B}(0,R)$. 
\end{lemma}
\begin{Proof}
    \noindent
$a)$ If $x \in B(0,R)$ we have 
$$
\mathbb{P}\bigg(\sup_{s\in [0,\tau^B(x))} |x+X_s| < R \bigg) = 1, 
$$
  as the distribution of the furthest reach from the origin immediately
  before exit time admits a density
  by Theorem~1.3-$(ii)$ of \cite{kyprianou2018}. 
  Similarly, letting
  $\widebar{B}^c(0,R):= \real^d \setminus \widebar{B}(0,R)$, 
  we have 
$$
\mathbb{P}\big(x+X_{\tau^B(x)} \in \widebar{B}^c(0,R)\big) = 1. 
$$
 as the distribution of
  $x+X_{\tau^B(x)}$ admits a density on $\real^d \setminus \widebar{B}(0,R)$ 
 by Equation~(2.2) in \cite{basscranston}. 
 Therefore, we have
 $$
 \P \bigg( \sup_{s\in [0,\tau^B(x))} |x+X_s| < R \quad
   \mbox{and}
   \quad
   x+X_{\tau^B(x)}\in \widebar{B}^c(0,R) \bigg) =1
    $$
   and almost surely
   there exists $r_0 (\omega )>0$ such that
$$
  ( y+X_s(\omega))_{ s \in [0,\tau^B(x))} 
  \subset B(0,R)  
    \quad
    \mbox{and}
    \quad 
y + X_{\tau^B(x)}(\omega) \in \widebar{B}^c(0,R), 
$$
  provided that $| y-x | < r_0 (\omega )$. 
  Therefore, for any $y \in B(x,r_0 (\omega ))$ we have
 $\tau^B(y)(\omega) = \tau^B(x)(\omega)$,  
 which proves \eqref{eq:tau}. 
\\
\noindent
$b)$
 If $x \in \widebar{B}(0,R) \setminus B(0,R)$ we have $\tau^B(x) = 0$,
 and by (6.3) in \cite{bogdanbarrier} there exists $C^*>0$
 such that for sufficiently small $\varepsilon$ we have 
$$
 \mathbb{P}\big(\tau^B(x_n) > \varepsilon \big) <  \frac{C^*}{\sqrt{\varepsilon}} V(d(x_n,B^c(0,R))),
  \qquad \varepsilon >0, 
$$
where $V$ is the renewal function of the ascending height-process,
which satisfies $\lim_{r\rightarrow 0} V(r)=0$, and
$B^c(0,R):= \real^d \setminus B(0,R)$. Therefore, if the
sequence $(x_n)_{n\geq 0}$ in $B(0,R)$ is such that 
\begin{eqnarray}
\nonumber 
V(d(x_n,B^c(0,R))) < \frac{1}{2^n}, \qquad n \geq 0, 
\end{eqnarray}
we have $ \sum_{n\geq 0} \mathbb{P}\big(\tau^B(x_n) > \varepsilon\big) < \infty$,
and we conclude by the Borel-Cantelli Lemma.
\end{Proof}
Based on \eqref{eq:tau}, we obtain a sufficient condition for the continuity of 
 $x\mapsto \mathbb{E}[\mathcal{H}(\mathcal{T}_{x})]$ 
in $ x\in \widebar{B}(0,R) \subset \real^d$.
We note that this continuity property may be extended to non spherical
domains for which \eqref{eq:tau} is satisfied. 
\begin{lemma}
\label{pl1}
 Assume  that   
  $(\mathcal{H}(\mathcal{T}_{x}))_{x\in B(0,R)}$ is uniformly integrable.
 Then, the function 
 $u(x):= \mathbb{E}[\mathcal{H}(\mathcal{T}_{x})]$ is
 continuous in $ x\in \widebar{B}(0,R)$.
 \end{lemma} 
\begin{Proof}
   Let $x\in \widebar{B}(0,R)$, and let $(x_n)_{n\geq 0}$ denote a sequence in
 $B(0,R)$ converging to $x$ and satisfying \eqref{eq:tau}.
 For any $\widebar{k} \in \mathcal{K}$
 we let $\tau_{{\widebar{k}},x} := \tau^B \big( {X^{\widebar{k}-}_{T_{\widebar{k}-},x}}\big)$,
 and note that the event  
$$
 A_{\widebar{k}} := 
 \left\{ \lim_{n \rightarrow \infty} \tau_{{\widebar{k}},x_n}  = \tau_{{\widebar{k}},x} \right\} \bigcap \left\{  \lim_{n \rightarrow \infty}X_{\cdot, x_n}^{\widebar{k}}= X_{\cdot, x}^{\widebar{k}} \right\}, 
$$
 has probability one
 by Lemma~\ref{l1} and
 the relation
 $X_{\cdot, x_n}^{\widebar{k}}= X_{\cdot, x}^{\widebar{k}} + x_n-x$,
 $n\geq 0$. 
By Lemma~\ref{l1}-$a)$, for some $n_0 (\omega )$ large enough we have 
$$
X_{\tau_{{\widebar{k}},x_n} }^{\widebar{k}}
=
X_{\tau_{{\widebar{k}},x} }^{\widebar{k}}
+x_n-x,  
$$
 and $\tau_{{\widebar{k}},x_n}=\tau_{{\widebar{k}},x}$, $n \geq n_0 (\omega )$. 
 Therefore, using the continuity of $\phi$ and $c_l,l\in \mathcal{L}$, we have 
$$
 \lim_{n \rightarrow \infty}  \phi\big(X_{\tau_{{\widebar{k}},x_n} }^{\widebar{k}}\big) \mathbbm{1}_{ \{ T_{\widebar{k}} = \tau_{{\widebar{k}},x_n} \} } =  \phi\big(X_{\tau_{{\widebar{k}},x}}^{\widebar{k}}\big) \mathbbm{1}_{\{ T_{\widebar{k}} = \tau_{{\widebar{k}},x}\}}, ~~ \mathbb{P}- \text{a.s.}
$$
 and
 $$
 \lim_{n \rightarrow \infty} \frac{c_{I_{\widebar{k}}}\big(X_{T_{\widebar{k}},x_n}^{\widebar{k}}\big)}{q_{I_{\widebar{k}}}} \mathbbm{1}_{ \{ T_{\widebar{k}} < \tau_{{\widebar{k}},x_n} \} } =  \frac{c_{I_{\widebar{k}}}\big(X_{T_{\widebar{k}},x}^{\widebar{k}}\big)}{q_{I_{\widebar{k}}}} \mathbbm{1}_{\{ T_{\widebar{k}} < \tau_{{\widebar{k}},x} \}}, ~~ \mathbb{P}- \text{a.s.}. 
$$
 Hence by \eqref{djsdiuoa}, 
 on the event $A := \bigcap_{\widebar{k} \in \mathcal{K} } A_{\widebar{k}}$ 
 of probability one, we have 
$$
\lim_{n \rightarrow \infty} \mathcal{H}(\mathcal{T}_{x_n}(\omega)) = \mathcal{H}(\mathcal{T}_{x}(\omega)). 
$$
Therefore,
 for any sequence $(x_n)_{n\geq 1}$ converging to
  $x\in \widebar{B}(0,R)$ fast enough, we have
$$
 \mathbb{P}\big(\lim_{n \rightarrow \infty} \mathcal{H}(\mathcal{T}_{x_n}) = \mathcal{H}(\mathcal{T}_{x})\big) = 1,
$$ 
 which yields $\lim_{n\to \infty} u(x_n) = u(x)$ by uniform integrability of
 $(\mathcal{H}(\mathcal{T}_{x}))_{x\in B(0,R)}$. 
\end{Proof}
\vspace{-0.2cm} 
  \begin{theorem}
   \label{t1.0}
 Assume that there exists a 
   sequence $(q_l)_{l \in \mathcal{L}}$ of positive numbers
   summing to one, such that the partial differential inequality 
   \begin{equation}
     \label{fjkdslf} 
    \Delta_s v (x) + \sum_{l \in \mathcal{L}} \frac{|c_l(x)|^p}{q_l^{p-1}} v^l (x) \leq v (x),
   \qquad x\in B(0,R),  
\end{equation}
admits a non-negative solution $v \in H^{2s} (B(0,R)) \bigcap \mathcal{C}\big(\widebar{B}(0,R)\big)$
such that $v \geq |\phi|^p$ on $\mathbb{R}^d \setminus B(0,R)$
for some $p>1$.  
 Then, the nonlinear PDE~\eqref{eq:1}
 admits a (continuous) viscosity solution on $\mathcal{O}=B(0,R)$. 
\end{theorem}
\begin{Proof} 
  We take $\rho(t) := e^{-t}, t\geq0$, so that $\mathcal{H}(\mathcal{T}_x)$
  rewrites as 
\begin{equation}
 \label{djsda} 
        {\mathcal{H} (\mathcal{T}_{x}) :=
         \prod_{\widebar{k} \in \mathcal{K}^{\circ}} \frac{c_{I_{\widebar{k}}}\big(X^{\widebar{k}}_{T_{\widebar{k}},x}\big)}{q_{I_{\widebar{k}}}} \prod_{\widebar{k} \in \mathcal{K}^{\partial}}   \phi\big(X^{\widebar{k}}_{T_{\widebar{k}},x}\big) }, \qquad x \in \real^d.
\end{equation}  
 Applying the It\^o-Dynkin formula and \eqref{fjkdslf}
 to the process $\big(e^{-t}v\big(X_{t,x}^{\widebar{1}}\big)\big)_{t\in \real_+}$
 on the time interval
 $[0,\tau^B(x)]$ we get 
$$
 v(x) \geq \mathbb{E} \left[ e^{-\tau^B(x)} \big|\phi\big(
   X_{\tau^B(x),x}^{\widebar{1}}
   \big)\big|^p + \int_0^{\tau^B(x)} e^{-t} \sum_{l \in \mathcal{L}} \frac{\big|c_l\big(
     X_{t,x}^{\widebar{1}}
     \big)\big|^p}{q_l^{p-1}}v^l\big(
   X_{t,x}^{\widebar{1}}
   \big)dt\right]. 
$$
 Thus, by the same recursion as in the proof of Proposition~\ref{p2}, we obtain 
$$
 v(x) \geq \mathbb{E} \left[ \prod_{\widebar{k} \in \bigcup_{i=1}^n \mathcal{K}_i^{\circ}} \frac{\big|c_{I_{\widebar{k}}}\big(X^{\widebar{k}}_{T_{\widebar{k}},x}\big)\big|^p}{q_{I_{\widebar{k}}}^p}  \prod_{\widebar{k} \in \bigcup_{i=1}^n  \mathcal{K}_i^{\partial}} \big|\phi\big(X^{\widebar{k}}_{T_{\widebar{k}},x} \big)\big|^p \prod_{\widebar{k} \in \mathcal{K}_{n+1}} v\big(X_{T_{\widebar{k}-}}^{\widebar{k}}\big) \right],
 \quad n \geq 1.
$$
 Letting $n$ tend to infinity and applying Fatou's Lemma,
 since $v$ is non-negative we find
$$
v(x) \geq \mathbb{E} \left[
  \prod_{\widebar{k} \in \mathcal{K}^\circ} \frac{\big|c_{I_{\widebar{k}}}\big(X^{\widebar{k}}_{T_{\widebar{k}},x}\big)\big|^p}{q_{I_{\widebar{k}}}^p}  \prod_{\widebar{k} \in \mathcal{K}^{\partial}} \big|\phi\big(X^{\widebar{k}}_{T_{\widebar{k}},x} \big)\big|^p \right]
= \mathbb{E} [ |\mathcal{H}(\mathcal{T}_{x}) |^p ],
\quad x \in B(0,R). 
$$
 In particular, $(\mathcal{H}(\mathcal{T}_{x}))_{x \in B(0,R)} $ is uniformly bounded
 in $L^p(B(0,R))$
 since $v \in \mathcal{C}\big(\widebar{B}(0,R)\big)$. 
 Therefore, by Proposition~\ref{p1} and Lemma~\ref{pl1}  
 $u$ is a continuous viscosity solution of the PDE~\eqref{eq:1} on $B(0,R)$. 
\end{Proof}
\noindent
To prove the integrability required in Theorems~\ref{t2}-\ref{t3} we adapt 
the approach of \cite{claisse} to the fractional setting, 
by constructing a branching process that stochastically dominates
the underlying stable branching process uniformly in $x\in B(0,R)$. 
\begin{Proofy}{\em\ref{t2}.}  
 As in the proof of Theorem~\ref{t1.0}, in order to show that the PDE~\eqref{eq:1}
 admits a (continuous) viscosity solution it suffices to show that 
 $(\mathcal{H}(\mathcal{T}_x))_{x\in B(0,R)}$ is uniformly bounded
 in $L^p(B(0,R))$ for some $p>1$. 
 We take again $\rho(t) := e^{-t}, t\geq0$, in which case $\mathcal{H}(\mathcal{T}_x)$
 rewrites as in \eqref{djsda}. Letting
 $$
 q_k : = \frac{
   |c_k|_{L^\infty (B(0,R))}
    }{\sum_{l \in \mathcal{L}} |c_l|_{L^\infty (B(0,R))}},
 \quad k \in \mathcal{L}, 
$$ 
 we have
 $| \mathcal{H} (\mathcal{T}_{x}) | \leq 1$
 provided that 
 $$
 C_0 := \max \left(
 |\phi|_{L^\infty (B^c(0,R))}, \sup_{l} \frac{|c_l|_{L^\infty (B(0,R))}}{q_l} \right)
 = \max\bigg(|\phi|_{L^\infty (B^c(0,R))} , \sum_{l \in \mathcal{L}} |c_l|_{L^\infty (B(0,R))} \bigg) \leq 1, 
$$
 in which case 
 $(\mathcal{H}(\mathcal{T}_{x}))_{x\in B(0,R)}$
 is uniformly integrable 
 and we conclude by
 Lemma~\ref{pl1} and Proposition~\ref{p1}.
 If $C_0>1$, we let 
\begin{equation} 
  \label{fjkds}
  \delta: = 1 - \inf_{x \in B(0,R)} \mathbb{E}\big[e^{-\tau^B(x)}\big], 
\end{equation}
 and
 $$
 \widetilde{f}(s) := \sum_{l \in \mathcal{L}} \widetilde{q}_l s^l
 $$
 where
 $$
 \widetilde{q}_0 := 1 - \delta + 
 \frac{\delta |c_0|_{L^\infty (B(0,R))}}{\sum_{l \in \mathcal{L}} |c_l|_{L^\infty (B(0,R))}},
 \quad \widetilde{q}_k := \frac{\delta |c_k|_{L^\infty (B(0,R))}}{\sum_{l \in \mathcal{L}} |c_l|_{L^\infty (B(0,R))}},
 \qquad
 k \geq 1.
 $$
 By Proposition~3.5 in \cite{claisse}, 
 $(\mathcal{H}(\mathcal{T}_x))_{x\in B(0,R)}$ is uniformly bounded in $L^p(B(0,R))$
 for some $p>1$ provided that $ C_0 < ( \gamma (s^*) )^{1/p}$, where 
$$
\gamma (s^*) :=
\frac{1}{\widetilde{f}'(s^*)}
= \frac{s^*}{\widetilde{f}(s^*)}
$$
and $s^*$ is the solution of 
$s^* \widetilde{f}'(s^*) = \widetilde{f}(s^*)$ if it exists,
or $s^* = \zeta$ otherwise,
where $\zeta$ denotes the radius of convergence of $\widetilde{f}$.
As above, we conclude the proof 
 from Lemma~\ref{pl1} and Proposition~\ref{p1}.
\end{Proofy}
\begin{Proofy}{\em\ref{t3}.}  
 When $\mathcal{L}$ is finite we have $\zeta = \infty$
and a solution $s^*$ to $s^* \widetilde{f}'(s^*) = \widetilde{f}(s^*)$ always
exists. 
By part~$(iv)$ of the proof of Proposition~3.5 in \cite{claisse},  
$s^*$ tends to infinity as $\delta$ goes to $0$
(or equivalently as ${\rm Diam} (B(0,R))$ goes to $0$ by \eqref{fjkds}) 
and $\lim_{\delta \to 0} \gamma (s^*) = \infty$.
Hence, if $R$ is sufficiently small
we have $C_0 < ( \gamma (s^*))^{1/p}$
and $(\mathcal{H}(\mathcal{T}_x))_{x\in B(0,R)}$ is uniformly bounded in $L^p(B(0,R))$
for some $p>1$ by Proposition~3.5 in \cite{claisse}.
 Therefore, we can conclude from 
 Lemma~\ref{pl1} and Proposition~\ref{p1}
 as in the proof of Theorem~\ref{t2}.
\end{Proofy}
\noindent
 In the next proposition we note that
 the viscosity solution 
 obtained in Theorems~\ref{t3} and \ref{t1.0} is also a weak solution,
 and that it is the only weak and viscosity solution
  of \eqref{eq:1} provided that $\phi$ belongs to $H^{2s} (\real^d)$. 
  \begin{prop}
  \label{p3}
  Assume that \eqref{eq:1} admits a viscosity solution $u$,
  and that $\phi$ belongs to $H^{2s} (\real^d)$. 
  Then $u \in H^s (\real^d)
  \cap \mathcal{C}\big(\widebar{B}(0,R)\big)$,
  it is a weak solution
  and the only weak and viscosity solution
  of \eqref{eq:1}. 
\end{prop} 
\begin{Proof}
  As $u$ is a viscosity solution it is continuous,
  and bounded on $\real^d$.
  Letting $v := u-\phi$, $v$ solves the equation 
 \begin{equation*}
\begin{cases}
  \displaystyle
  \Delta_s v(x) + \Delta_s \phi(x) +
 \sum\limits_{l \in \mathcal{L}} c_l(x) (v+\phi)^{l}(x)
  = v(x)+\phi(x), \quad x \in B(0,R), 
  \medskip
  \\
v(x) = 0, \quad x \in \mathbb{R}^d \setminus B(0,R),
\end{cases}
\end{equation*}
which can be rewritten as:
 \begin{equation}
 \label{eq:2}
\begin{cases}
  \displaystyle
  \Delta_s v(x) +
 \sum\limits_{l \in \mathcal{L}} \widetilde{c}_l(x) v^{l}(x)
  = 0, \quad x \in B(0,R), 
  \medskip
  \\
v(x) = 0, \quad x \in \mathbb{R}^d \setminus B(0,R),
\end{cases}
\end{equation}
 where $\widetilde{c}_l$
 is a continuous bounded function 
 on $\real^d$ for every $l\in \mathcal{L}$. 
 Letting $g(x) := \sum\limits_{l \in \mathcal{L}} \widetilde{c}_l(x) v^{l}(x) $, 
 the Dirichlet problem
 \begin{equation*}
\begin{cases}
  \displaystyle
 \Delta_s w
  = - g(x), \quad x \in B(0,R), 
  \medskip
  \\
w(x) = 0, \quad x \in \mathbb{R}^d \setminus B(0,R),
\end{cases}
\end{equation*}
 admits $v$ as viscosity solution.
 Since $g$ is bounded on $B(0,R)$,
 it also admits a unique weak solution $w$
 which is a viscosity solution by Remark~2.11 in \cite{ros-oton}.
 By Theorem~5.2 in \cite{caffarelli},
 the viscosity solutions $v,w$ coincide
 and $v$ is the unique solution of \eqref{eq:2} in
 both in the weak and viscosity senses. 
 Therefore,
 $u$ is the unique solution of \eqref{eq:1}.
 In addition,
 since $u$ is a weak solution it belongs to $H^s(\real^d)$,
 see Proposition~1.4 in \cite{ros-oton}.
\end{Proof}

\tikzstyle{level 1}=[level distance=4cm, sibling distance=3cm]
\tikzstyle{level 2}=[level distance=5cm, sibling distance=5cm]

\section{Numerical examples} 
\label{s5}
   In this section we consider numerical examples involving the fractional Laplacian $\Delta_s$ and
 the $s$-stable subordinator
$(S_t)_{t\in \real_+}$
 with Laplace exponent $\eta(\lambda) = (2\lambda )^s$
 for $s \in (1/2,1)$.
 We represent $(X_t)_{t\in \real_+}$ 
 using the subordination $X_t := B_{S_t}$, where
 $(B_t)_{t\in \real_+}$ is a standard $d$-dimensional Brownian motion
 and $(S_t)_{t\in \real_+}$
 is a L\'evy subordinator with Laplace exponent $\eta$,
 defined by 
$$
 \mathbb{E}\big[e^{-\lambda S_t}\big] = e^{-t ( 2 \lambda)^s}, \qquad
\lambda , \ \! t \geq 0, 
$$
 see e.g. Theorem~1.3.23 in \cite{applebk2}. 
 For the generation of random samples of $S_t$, we use the formula
 
$$
\widetilde{S}_t := 2t^{1/s}\frac{\sin( s \big(U+\pi / 2) \big)}{\cos^{1/s} (U)} \left(
\frac{\cos\big(U- s  (U+\pi / 2) \big)}{E}\right)^{-1+1/s }
$$
based on the
Chambers-Mallows-Stuck (CMS) method, 
where $U \sim U (-\pi / 2,\pi /2)$,
and $E \sim {\rm Exp} (1)$, 
see Relation~(3.2) in \cite{weron1996chambers}, 
where $\psi(\lambda)$
denotes the L\'evy symbol of $(S_t)_{t\in \real_+}$. 
For $k\geq 0$, we consider the function 
$$
\Phi_{k,s } (x) := (1-| x|^2)^{k+s}_+,
\qquad x\in \real^d, 
$$
 which is Lipschitz if $k>1-s$, and solves the Poisson problem  
 $\Delta_s \Phi_{k,s} = -\Psi_{k,s}$
 on $\real^d$, with 
 \begin{align} 
\nonumber 
&  \Psi_{k,s} (x)
\\
\nonumber 
& := \left\{
\begin{array}{ll}
  \displaystyle
 \frac{\Gamma(s+d/2)
   \Gamma(k+1+ s )}{
   4^{-s}
   \Gamma(k+1)\Gamma( d / 2 )} 
~{}_2F_1\left(
  \frac{d}{2} + s ,-k;\frac{d}{2};| x|^2
  \right), ~~ | x| \leq 1
  \medskip 
  \\
\displaystyle
 \frac{4^s \Gamma( s+d/2 )
  \Gamma(k+1+ s )}{\Gamma(k+1+ s+ d/2 )
   \Gamma(- s )
   | x|^{d+2s}
 } 
  {}_2F_1\left(
  \frac{d}{2} + s , 1 + s ;k+1+\frac{d}{2} + s ;
  \frac{1}{| x|^2}
\right), ~~| x|>1
\end{array}
\right.
\end{align} 
 $x\in \real^d$, where ${}_2F_1 ( a,b;c;y)$ is 
 Gauss's hypergeometric function,
 see (5.2) in \cite{getoor} and 
 Lemma~4.1 in \cite{biler2015nonlocal}, 
 and Relation~(36) in \cite{oberman}. 

\subsubsection*{Dirichlet problem} 
 We solve the Dirichlet problem
\begin{equation}\label{eq:nld0}
\begin{cases}
  \displaystyle
 \Delta_s u (x) + \Psi_{k,s}(x) = 0, ~~ x \in B(0,1)
  \medskip
  \\
u(x) =0, ~~ x \in B^c(0,1). 
\end{cases}{}
\end{equation}{}
 with 
 $$
 c_0(x) u^0(x):= \Psi_{k,s}(x),
 \qquad
 c_1(x) u^1(x) := u(x), 
$$
 and explicit solution 
\begin{equation}
 \label{djk1}
   u(x) = \Phi_{k,s} (x) 
   = (1-| x|^2)^{k+s}_+, \qquad
x \in \real^d. 
\end{equation} 
We note that when $k=0$, \eqref{eq:nld0} reads
$$
\Delta_s u (x) +
\frac{\Gamma(s+d/2)
   \Gamma(s +1)}{
   4^{-s}
   \Gamma( d / 2 )} 
=0,
$$
while when $k=1$ we have
$$
\Delta_s u (x) +
 \frac{\Gamma(s+d/2)
   \Gamma(s +2)}{
   4^{-s}
   \Gamma(2)\Gamma( d / 2 )} 
 \left(
  1 - \left( 1+ \frac{2s}{d} \right) | x|^2
  \right) = 0.
  $$
 The random tree associated to \eqref{eq:nld0}
starts at a point $x\in B(0,1)$
and branches into \textcolor{cyan}{zero branch} or \textcolor{blue}{one branch} 
 as in the following random samples: 
\bigskip 
\begin{center}
\resizebox{0.90\textwidth}{!}{
\begin{tikzpicture}[scale=0.9,grow=right, sloped]
\node[ellipse split,draw,blue,text=black,thick]{$0$ \nodepart{lower} $x$}
    child {
        node[ellipse split,draw,blue,text=black,thick] {$T_{\bar{1}}$ \nodepart{lower} $X^{\bar{1}}_{T_{\bar{1}},x}$}  
            child {
                node[ellipse split,draw,blue,text=black,thick, right=0.cm] {$T_{(1,1)}$ \nodepart{lower} $X^{(1,1)}_{T_{(1,1)},x}$}
                child{
                    node[ellipse split,draw,black, right=0.cm] {$T_{(1,1,1)}:=T_{(1,1)}+\tau^B(X^{(1,1,1)}_{T_{(1,1,1)},x})$ \nodepart{lower} $X^{(1,1,1)}_{T_{(1,1,1)},x}$}
                    edge from parent
                    node[above]{~$(1,1,1)$~}
            }
                edge from parent
                node[above] {$(1,1)$}
            }
            edge from parent
            node[above] {$\bar{1}$}
    };
\end{tikzpicture}
}
\end{center}
\begin{center}
\resizebox{0.5\textwidth}{!}{
\begin{tikzpicture}[scale=0.9,grow=right, sloped]
\node[ellipse split,draw,blue,text=black,thick]{$0$ \nodepart{lower} $x$}
    child {
        node[ellipse split,draw,blue,text=black,thick] {$T_{\bar{1}}$ \nodepart{lower} $X^{\bar{1}}_{T_{\bar{1}},x}$}  
            child {
                node[ellipse split,draw,cyan,text=black,thick, right=0.cm] {$T_{(1,1)}$ \nodepart{lower} $X^{(1,1)}_{T_{(1,1)},x}$
            }
                edge from parent
                node[above] {$(1,1)$}
            }
            edge from parent
            node[above] {$\bar{1}$}
    };
\end{tikzpicture}
}
\end{center}

\begin{figure}[H]
\centering
\hskip-0.4cm
\begin{subfigure}{.5\textwidth}
\centering
\includegraphics[width=\textwidth]{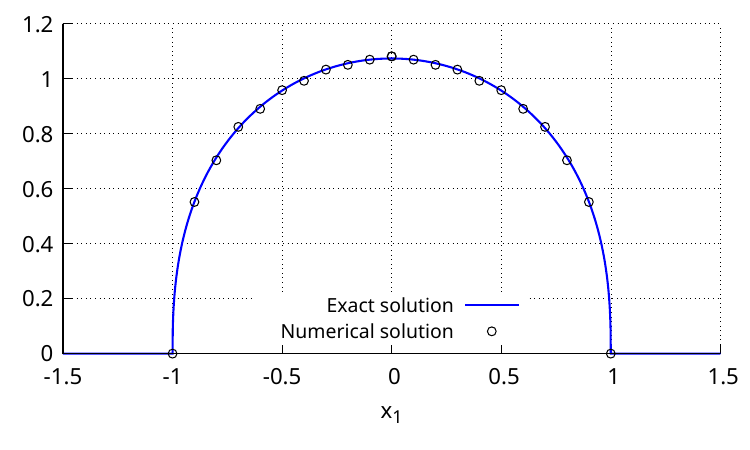}
\vskip-0.3cm
\caption{Numerical solution of
  $\Delta_s u (x) = -1$ with~$k$=0.}
\end{subfigure}
\begin{subfigure}{.50\textwidth}
\centering
\includegraphics[width=\textwidth]{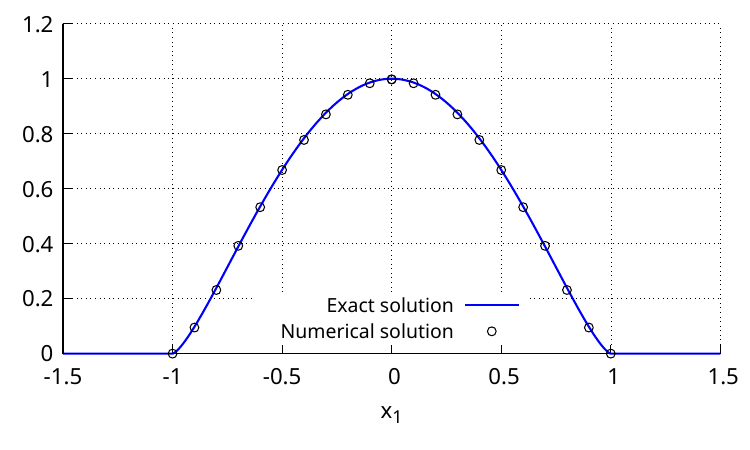}
\vskip-0.3cm
\caption{Numerical solution of \eqref{eq:nld0} with $k=1$.}
\end{subfigure}
\caption{Numerical solutions 
   in dimension $d=1$ with $\alpha =0.8$.}
\label{fig1}
\end{figure}
\vskip-0.3cm
\noindent 
Figure~\ref{fig1}-$(a)$, which uses one million Monte Carlo samples,
 can be compared to Figure~6.5a in \cite{huang-oberman}.

\begin{figure}[H]
\centering
\hskip-0.4cm
\begin{subfigure}{.5\textwidth}
\centering
\includegraphics[width=\textwidth]{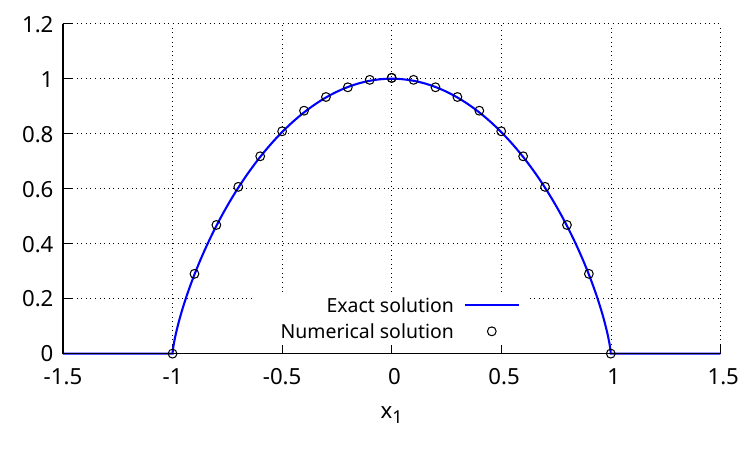}
\vskip-0.3cm
\caption{Numerical solution of \eqref{eq:nld0} with $k=0$.}
\end{subfigure}
\begin{subfigure}{.50\textwidth}
\centering
\includegraphics[width=\textwidth]{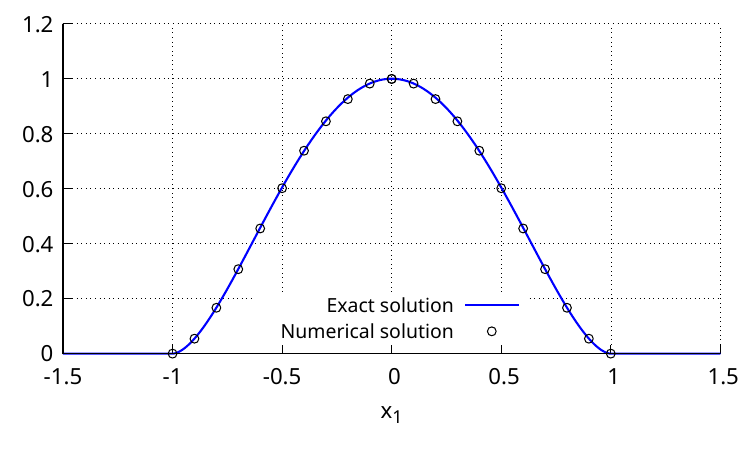}
\vskip-0.3cm
\caption{Numerical solution of \eqref{eq:nld0} with $k=1$.}
\end{subfigure}
\caption{Numerical solutions of \eqref{eq:nld0} in dimension $d=10$ with $\alpha =1.5$.}
\label{fig1-0}
\end{figure}
\vskip-0.3cm
\noindent 
\subsubsection*{Linear fractional elliptic PDE}
\begin{figure}[H]
\centering
\hskip-0.4cm
\begin{subfigure}{.5\textwidth}
\centering
\includegraphics[width=\textwidth]{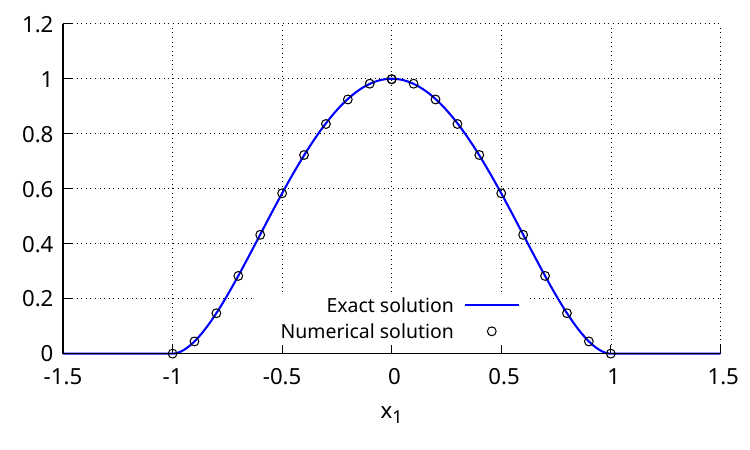}
\vskip-0.3cm
\caption{Numerical solution of \eqref{eq:nld1} with $k=1$.}
\end{subfigure}
\begin{subfigure}{.50\textwidth}
\centering
\includegraphics[width=\textwidth]{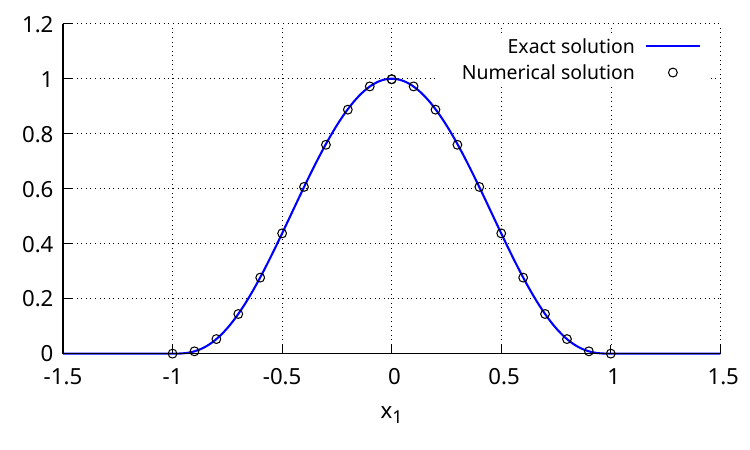}
\vskip-0.3cm
\caption{Numerical solution of \eqref{eq:nld1} with $k=2$.}
\end{subfigure}
\caption{Numerical solutions of \eqref{eq:nld1} in dimension $d=100$ with $\alpha =1.75$.}
\label{fig2}
\end{figure}

 We solve the linear elliptic problem 
\begin{equation}\label{eq:nld1}
\begin{cases}
  \displaystyle
 \Delta_s (x) + \Psi_{k,s}(x) -  (1-| x|^2)^{k+s}_+ + u (x) =0, ~~ x \in B(0,1)
  \medskip
  \\
u(x) =0, ~~ x \in B^c(0,1), 
\end{cases}{}
\end{equation}{}
 with 
$$ 
 c_0(x) u^0(x):= \Psi_{k,s}(x) -  (1-| x|^2)^{2k+2s}_+,
 \qquad
 c_1(x) u^1(x) := 2u(x),
 $$
 and explicit solution \eqref{djk1}. 
 The random tree associated to \eqref{eq:nld1}
 is the same as in the previous example, 
 and the simulations of Figure~\ref{fig2}
 use five million Monte Carlo samples. 

\subsubsection*{Nonlinear fractional elliptic PDE}
\begin{figure}[H]
\centering
\hskip-0.4cm
\begin{subfigure}{.5\textwidth}
\centering
\includegraphics[width=\textwidth]{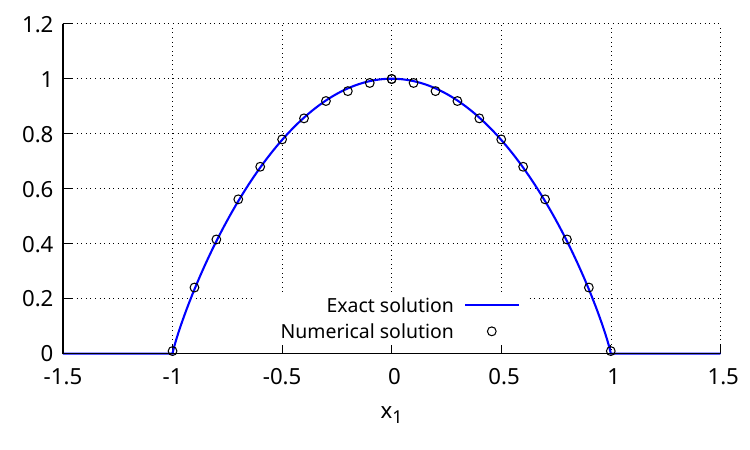}
\vskip-0.3cm
\caption{Numerical solution of \eqref{eq:nld2} with $k=0$.}
\end{subfigure}
\begin{subfigure}{.50\textwidth}
\centering
\includegraphics[width=\textwidth]{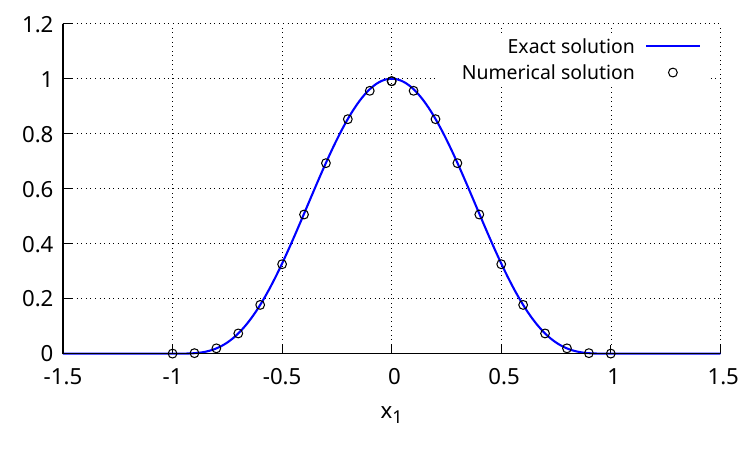}
\vskip-0.3cm
\caption{Numerical solution of \eqref{eq:nld2} with $k=3$.}
\end{subfigure}
\caption{Numerical solutions of \eqref{eq:nld2} in dimension $d=10$ with $\alpha =1.75$.}
\label{fig3}
\end{figure}

 Here we aim at recovering the explicit solution
\eqref{djk1}
 of the nonlinear elliptic PDE
\begin{equation}\label{eq:nld2}
\begin{cases}
\displaystyle
 \Delta_s u (x) + \Psi_{k,s}(x) -  (1-| x|^2)^{2k+2s}_+ + u^2 (x) =0, ~~ x \in B(0,1),
  \medskip
  \\
u(x) =0, ~~ x \in B^c(0,1), 
\end{cases}{}
\end{equation}{}
\noindent
 with 
$$ 
 c_0(x) u^0(x)= \Psi_{k,s}(x) -  (1-| x|^2)^{2k+2s}_+,
 \quad 
 c_1(x) u^1(x) = u(x),
 \quad 
c_2(x) u^2(x) = u^2(x). 
$$
 The random tree associated to \eqref{eq:nld2} starts at point $x\in B(0,1)$ and branches into \textcolor{cyan}{zero branch}, \textcolor{blue}{one branch},
 or \textcolor{violet}{two branches}, as in the following random sample tree: 
 
 \medskip

\begin{center}
\resizebox{0.85\textwidth}{!}{
\begin{tikzpicture}[scale=0.9,grow=right, sloped]
\node[ellipse split,draw,blue,text=black,thick]{$0$ \nodepart{lower} $x$}
    child {
        node[ellipse split,draw,violet,text=black,thick] {$T_{\bar{1}}$ \nodepart{lower} $X^{\bar{1}}_{T_{\bar{1}},x}$}   
            child {
                node[ellipse split,draw,violet,text=black,thick] {$T_{(1,2)}$ \nodepart{lower} $X^{(1,2)}_{T_{(1,2)},x}$} 
                child{
                node[ellipse split,draw,cyan,text=black,thick,right=5cm,below=-2.1cm]{$T_{(1,2,2)}$ \nodepart{lower} $X^{(1,2,2)}_{T_{(1,2,2)},x}$}
                edge from parent
                node[above]{$(1,2,2)$}
                }
                child{
                node[ellipse split,draw,thick, right=4.97cm, below=-0.6cm]{$T_{(1,2,1)}:= T_{(1,2)}+ \tau^B(X^{(1,2,1)}_{T_{(1,2,1)},x})$ \nodepart{lower} $X^{(1,2,1)}_{T_{(1,2,1)},x} $}
                edge from parent
                node[above]{$(1,2,1)$}
                }
                edge from parent
                node[above] {$(1,2)$}
            }
            child {
                node[ellipse split,draw,blue,text=black,thick, right=0.cm] {$T_{(1,1)}$ \nodepart{lower} $X^{(1,1)}_{T_{(1,1)},x}$}
                child{
                    node[ellipse split,draw,black, right=0.cm] {$T_{(1,1,1)} := T_{(1,1)}+\tau^B(X^{(1,1,1)}_{T_{(1,1,1)},x})$ \nodepart{lower} $X^{(1,1,1)}_{T_{(1,1,1)},x}$}
                    edge from parent
                    node[above]{~$(1,1,1)$~}
            }
                edge from parent
                node[above] {$(1,1)$}
            }
            edge from parent 
            node[above] {$\bar{1}$}
    };
\end{tikzpicture}
} 
\end{center}
The simulations of Figure~\ref{fig3} use 20 million Monte Carlo samples.

\footnotesize

\def\cprime{$'$} \def\polhk#1{\setbox0=\hbox{#1}{\ooalign{\hidewidth
  \lower1.5ex\hbox{`}\hidewidth\crcr\unhbox0}}}
  \def\polhk#1{\setbox0=\hbox{#1}{\ooalign{\hidewidth
  \lower1.5ex\hbox{`}\hidewidth\crcr\unhbox0}}} \def\cprime{$'$}


\begin{thebibliography}{35}
\providecommand{\natexlab}[1]{#1}
\providecommand{\url}[1]{\texttt{#1}}
\expandafter\ifx\csname urlstyle\endcsname\relax
  \providecommand{\doi}[1]{doi: #1}\else
  \providecommand{\doi}{doi: \begingroup \urlstyle{rm}\Url}\fi

\bibitem[Agarwal and Claisse(2020)]{claisse}
A.~Agarwal and J.~Claisse.
\newblock Branching diffusion representation of semi-linear elliptic {PDE}s and
  estimation using {M}onte {C}arlo method.
\newblock \emph{Stochastic Processes and their Applications}, 130\penalty0
  (8):\penalty0 5006--5036, 2020.

\bibitem[Alves et~al.(2020)Alves, Bisci, and Ledesma]{alves}
C.O. Alves, G.~M. Bisci, and C.E.~Torres Ledesma.
\newblock Existence of solutions for a class of fractional problems on exterior
  domains.
\newblock \emph{J. Differential Equations}, 268\penalty0 (11):\penalty0
  7183--7219, 2020.

\bibitem[Applebaum(2009)]{applebk2}
D.~Applebaum.
\newblock \emph{L\'{e}vy processes and stochastic calculus}, volume 116 of
  \emph{Cambridge Studies in Advanced Mathematics}.
\newblock Cambridge University Press, Cambridge, second edition, 2009.

\bibitem[Barles et~al.(2008)Barles, Chasseigne, and Imbert]{barles2}
G.~Barles, E.~Chasseigne, and C.~Imbert.
\newblock On the {D}irichlet problem for second-order elliptic
  integro-differential equations.
\newblock \emph{Indiana Univ. Math. J.}, 57\penalty0 (1):\penalty0 213--246,
  2008.

\bibitem[Bass and Cranston(1983)]{basscranston}
R.F. Bass and M.~Cranston.
\newblock Exit times for symmetric stable processes in {${\bf R}^{n}$}.
\newblock \emph{Ann. Probab.}, 11\penalty0 (3):\penalty0 578--588, 1983.

\bibitem[Belak et~al.(2020)Belak, Hoffmann, and Seifried]{belak}
C.~Belak, D.~Hoffmann, and F.T. Seifried.
\newblock Branching diffusions with jumps and valuation with systemic
  counterparties.
\newblock Available at SSRN: https://ssrn.com/abstract=3451280 or
  https://dx.doi.org/10.2139/ssrn.3451280, 2020.

\bibitem[Biler et~al.(2015)Biler, Imbert, and Karch]{biler2015nonlocal}
P.~Biler, C.~Imbert, and G.~Karch.
\newblock The nonlocal porous medium equation: Barenblatt profiles and other
  weak solutions.
\newblock \emph{Archive for Rational Mechanics and Analysis}, 215\penalty0
  (2):\penalty0 497--529, 2015.

\bibitem[Bogdan et~al.(2015)Bogdan, Grzywny, and Ryznar]{bogdanbarrier}
K.~Bogdan, T.~Grzywny, and M.~Ryznar.
\newblock Barriers, exit time and survival probability for unimodal {L}\'{e}vy
  processes.
\newblock \emph{Probab. Theory Related Fields}, 162\penalty0 (1-2):\penalty0
  155--198, 2015.

\bibitem[Bony et~al.(1968)Bony, Courr\`ege, and Priouret]{bony}
J.-M. Bony, P.~Courr\`ege, and P.~Priouret.
\newblock Semi-groupes de {F}eller sur une vari\'{e}t\'{e} \`a bord compacte et
  probl\`emes aux limites int\'{e}gro-diff\'{e}rentiels du second ordre donnant
  lieu au principe du maximum.
\newblock \emph{Ann. Inst. Fourier (Grenoble)}, 18\penalty0 (fasc. 2):\penalty0
  369--521 (1969), 1968.

\bibitem[Caffarelli and Silvestre(2009)]{caffarelli}
L.~Caffarelli and L.~Silvestre.
\newblock Regularity theory for fully nonlinear integro-differential equations.
\newblock \emph{Comm. Pure Appl. Math.}, 62\penalty0 (5):\penalty0 597--638,
  2009.

\bibitem[Correia and Oliveira(2022)]{correia}
J.N. Correia and C.P. Oliveira.
\newblock Existence of a positive solution for a class of fractional elliptic
  problems in exterior domains involving critical growth.
\newblock \emph{J. Math. Anal. Appl.}, 506, 2022.

\bibitem[de~Pablo and S{\'a}nchez(2010)]{pablo}
A.~de~Pablo and U.~S{\'a}nchez.
\newblock Some {L}iouville-type results for a fractional equation.
\newblock Preprint, 2010.

\bibitem[Fall and Weth(2012)]{fall}
M.M. Fall and T.~Weth.
\newblock Nonexistence results for a class of fractional elliptic boundary
  value problems.
\newblock \emph{J. Funct. Anal.}, 263\penalty0 (8):\penalty0 2205--2227, 2012.

\bibitem[Felsinger et~al.(2015)Felsinger, Kassmann, and Voigt]{felsinger}
M.~Felsinger, M.~Kassmann, and P.~Voigt.
\newblock The {D}irichlet problem for nonlocal operators.
\newblock \emph{Math. Z.}, 279\penalty0 (3-4):\penalty0 779--809, 2015.

\bibitem[Getoor(1961)]{getoor}
R.K. Getoor.
\newblock First passage times for symmetric stable processes in space.
\newblock \emph{Trans. Amer. Math. Soc.}, 101:\penalty0 75--90, 1961.

\bibitem[Henry-Labord\`ere et~al.(2019)Henry-Labord\`ere, Oudjane, Tan, Touzi,
  and Warin]{labordere}
P.~Henry-Labord\`ere, N.~Oudjane, X.~Tan, N.~Touzi, and X.~Warin.
\newblock Branching diffusion representation of semilinear {PDE}s and {M}onte
  {C}arlo approximation.
\newblock \emph{Ann. Inst. H. Poincar\'e Probab. Statist.}, 55\penalty0
  (1):\penalty0 184--210, 2019.

\bibitem[Huang and Oberman(2014)]{huang-oberman}
Y.~Huang and A.~Oberman.
\newblock Numerical methods for the fractional {L}aplacian: a finite
  difference-quadrature approach.
\newblock \emph{SIAM J. Numer. Anal.}, 52\penalty0 (6):\penalty0 3056--3084,
  2014.

\bibitem[Huang and Oberman(2016)]{oberman}
Y.~Huang and A.~Oberman.
\newblock Finite difference methods for fractional {L}aplacians.
\newblock Preprint arXiv:1611.00164, 2016.

\bibitem[Ikeda et~al.(1968-1969)Ikeda, Nagasawa, and Watanabe]{inw}
N.~Ikeda, M.~Nagasawa, and S.~Watanabe.
\newblock Branching {M}arkov processes {I}, {II}, {III}.
\newblock \emph{J. Math. Kyoto Univ.}, 8-9:\penalty0 233--278, 365--410,
  95--160, 1968-1969.

\bibitem[Kwa{\' s}nicki(2017)]{tendef}
M.~Kwa{\' s}nicki.
\newblock Ten equivalent definitions of the fractional {L}aplace operator.
\newblock \emph{Fractional Calculus and Applied Analysis}, 20\penalty0
  (1):\penalty0 7--51, 2017.

\bibitem[Kyprianou et~al.(2020)Kyprianou, Rivero, and
  Satitkanitkul]{kyprianou2018}
A.E. Kyprianou, V.~Rivero, and W.~Satitkanitkul.
\newblock Deep factorisation of the stable process {III}: {T}he view from
  radial excursion theory and the point of closest reach.
\newblock \emph{Potential Anal.}, 53\penalty0 (4):\penalty0 1347--1375, 2020.

\bibitem[{Le Gall}(1995)]{LGBroSna}
J.-F. {Le Gall}.
\newblock The {B}rownian snake and solutions of {$\Delta u=u^2$} in a domain.
\newblock \emph{Probab. Theory Related Fields}, 102\penalty0 (3):\penalty0
  393--432, 1995.

\bibitem[L{\'o}pez-Mimbela(1996)]{lm}
J.A. L{\'o}pez-Mimbela.
\newblock A probabilistic approach to existence of global solutions of a system
  of nonlinear differential equations.
\newblock In \emph{Fourth Symposium on Probability Theory and Stochastic
  Processes (Spanish) (Guanajuato, 1996)}, volume~12 of \emph{Aportaciones Mat.
  Notas Investigaci\'on}, pages 147--155. Soc. Mat. Mexicana, M\'exico, 1996.

\bibitem[Mou(2017)]{mou}
C.~Mou.
\newblock Perron's method for nonlocal fully nonlinear equations.
\newblock \emph{Analysis and PDE}, 10\penalty0 (5):\penalty0 1227--1254, 2017.

\bibitem[Nagasawa and Sirao(1969)]{N-S}
M.~Nagasawa and T.~Sirao.
\newblock Probabilistic treatment of the blowing up of solutions for a
  nonlinear integral equation.
\newblock \emph{Trans. Amer. Math. Soc.}, 139:\penalty0 301--310, 1969.

\bibitem[Penent and Privault(2022)]{penent}
G.~Penent and N.~Privault.
\newblock Existence and probabilistic representation of the solutions of
  semilinear parabolic {PDE}s with fractional {L}aplacians.
\newblock \emph{Stochastics and Partial Differential Equations: Analysis and
  Computations}, 10:\penalty0 446--474, 2022.

\bibitem[Ros-Oton(2016)]{ros-oton2016}
X.~Ros-Oton.
\newblock Nonlocal elliptic equations in bounded domains: a survey.
\newblock \emph{Publ. Mat.}, 60\penalty0 (1):\penalty0 3--26, 2016.

\bibitem[Ros-Oton and Serra(2012)]{ro2012pohone}
X.~Ros-Oton and J.~Serra.
\newblock Fractional {L}aplacian: {P}ohozaev identity and nonexistence results.
\newblock \emph{C. R. Math. Acad. Sci. Paris}, 350\penalty0 (9-10):\penalty0
  505--508, 2012.

\bibitem[Ros-Oton and Serra(2014{\natexlab{a}})]{ro2014poho}
X.~Ros-Oton and J.~Serra.
\newblock The {P}ohozaev identity for the fractional {L}aplacian.
\newblock \emph{Arch. Ration. Mech. Anal.}, 213\penalty0 (2):\penalty0
  587--628, 2014{\natexlab{a}}.

\bibitem[Ros-Oton and Serra(2014{\natexlab{b}})]{ros-oton}
X.~Ros-Oton and J.~Serra.
\newblock The {D}irichlet problem for the fractional {L}aplacian: Regularity up
  to the boundary.
\newblock \emph{J. Math. Pures Appl.}, 101\penalty0 (3):\penalty0 275--302,
  2014{\natexlab{b}}.

\bibitem[Servadei and Valdinoci(2012)]{servadei}
R.~Servadei and E.~Valdinoci.
\newblock Mountain pass solutions for non-local elliptic operators.
\newblock \emph{J. Math. Anal. Appl.}, 389\penalty0 (2):\penalty0 887--898,
  2012.

\bibitem[Servadei and Valdinoci(2014)]{servadei2014}
R.~Servadei and E.~Valdinoci.
\newblock Weak and viscosity solutions of the fractional {L}aplace equation.
\newblock \emph{Publ. Mat.}, 58\penalty0 (1):\penalty0 133--154, 2014.

\bibitem[Skorokhod(1964)]{skorohodbranching}
A.V. Skorokhod.
\newblock Branching diffusion processes.
\newblock \emph{Teor. Verojatnost. i. Primenen.}, 9:\penalty0 492--497, 1964.

\bibitem[Stegli{\'n}ski(2021)]{steglinski}
R.~Stegli{\'n}ski.
\newblock Existence of a unique solution to a fractional partial differential
  equation and its continuous dependence on parameters.
\newblock \emph{Entropy}, 23, 2021.

\bibitem[Weron(1996)]{weron1996chambers}
R.~Weron.
\newblock On the {C}hambers-{M}allows-{S}tuck method for simulating skewed
  stable random variables.
\newblock \emph{Statistics \& probability letters}, 28\penalty0 (2):\penalty0
  165--171, 1996.

\end{thebibliography}
\end{document}